\documentclass[12pt]{amsart}
\usepackage{amssymb,amsthm,rotating, mathtools}
\usepackage[all]{xy}

\numberwithin{equation}{section}
 \newtheorem{thm}[equation]{Theorem}
\newtheorem{theorem}[equation]{Theorem}
 \newtheorem{defn}[equation]{Definition}
 \newtheorem{prop}[equation]{Proposition}

\newtheorem{cor}[equation]{Corollary}
 \newtheorem{lemma}[equation]{Lemma}

 \theoremstyle{definition}

\newtheorem{remark}[equation]{Remark}

 \newtheorem{example}[equation]{Example}

\DeclareMathOperator{\Ext}{Ext}
\DeclareMathOperator{\Hom}{Hom}
\DeclareMathOperator{\codim}{codim}
\DeclareMathOperator{\GL}{GL}
\DeclareMathOperator{\Gl}{GL}

\DeclareMathOperator{\Ima}{im}
\DeclareMathOperator{\Proj}{Proj}
\DeclareMathOperator{\proj}{proj}
\DeclareMathOperator{\sgn}{sgn}
\DeclareMathOperator{\Span}{Span}
\DeclareMathOperator{\Sym}{Sym} 
\DeclareMathOperator{\vol}{vol}
\DeclareMathOperator{\diag}{diag}
\DeclareMathOperator{\re}{res}

\newcommand{\HHD}{{\rm HH}^{\DOT}}
\newcommand{\rey}{\mathcal{R}}
\newcommand{\res}{{\re}}
\newcommand{\resAS}{{\re}}

\newcommand{\projH}{{\Proj}_H} 
\newcommand{\DOT}{\setlength{\unitlength}{1pt}\begin{picture}(2.5,2)
                  (1,1)\put(2,3.5){\circle*{2}}\end{picture}}

\newcommand{\del}{\partial}
\newcommand{\HH}{{\rm HH}}
\newcommand{\Wedge}{\textstyle\bigwedge}
\newcommand{\CC}{\mathbb{C}}

\newcommand{\ot}{\otimes}

\newcommand{\cohg}{{\rm HH}^{\DOT}\bigl(S(V),S(V)\bar{g}\bigr)}
\newcommand{\ep}{\epsilon}

\newcommand{\scooch}{\hspace{-.25ex}}
\newcommand{\ignore}[1]{\relax}

\newcommand{\rh}{\Theta} 
\newcommand{\ta}{\Upsilon\hspace{-.35ex}} 
\newcommand{\gamm}{\Gamma} 
\newcommand{\ttau}{{\Upsilon}\hspace{-1.3ex}\rule[.7ex]{.8ex}{.1ex}\hspace{.1ex}}

\title[Group actions on algebras]
{Group actions on algebras and the graded Lie structure of
Hochschild cohomology}

\date{November 19, 2010}


\author{Anne V.\ Shepler}
\address{Department of Mathematics, University of North Texas,
Denton, Texas 76203, USA}
\email{ashepler@unt.edu}
\author{Sarah Witherspoon}
\address{Department of Mathematics\\Texas A\&M University\\
College Station, Texas 77843, USA}\email{sjw@math.tamu.edu}
\thanks{The first author was partially supported by NSF grants
\#DMS-0402819 and \#DMS-0800951 and a research fellowship from the 
Alexander von Humboldt Foundation. 
The second author was partially supported by NSA grant
\#H98230-07-1-0038 and NSF grant
\#DMS-0800832.
Both authors were jointly supported by Advanced Research Program Grant 010366-0046-2007
from the Texas Higher Education Coordinating Board.}

\begin{document}
\maketitle

\begin{abstract}
Hochschild cohomology  governs deformations of algebras, and its graded 
Lie structure plays a vital role.
We study this structure for the Hochschild cohomology of the skew group algebra
formed by a finite group acting on an algebra by automorphisms.
We examine the Gerstenhaber bracket with a view toward deformations
 and developing  bracket formulas. 
We then focus on the linear group actions and polynomial algebras
that arise in orbifold theory and representation theory; 
deformations in this context include graded Hecke algebras and
symplectic reflection algebras.
We  give some general results describing when brackets are zero for
polynomial skew group algebras, which allow us in particular to find 
noncommutative Poisson structures.
For abelian groups, we express the bracket 
using inner products of group characters. Lastly,
we interpret results for graded Hecke algebras.
\end{abstract}

%%%%%%%%%%%%%%%%%%%%%%%%%%%%%%%%%%%%%%%%%%%%%%%%%%%%%%%%%%%%%%
%%%%%%%%%%%%%%%%%%%%%%%%%%%% SECTION %%%%%%%%%%%%%%%%%%%%%%%%%
%%%%%%%%%%%%%%%%%%%%%%%%%%%%%%%%%%%%%%%%%%%%%%%%%%%%%%%%%%%%%%
\section{Introduction}

Let $G$ be a finite group acting on a $\CC$-algebra $S$ by automorphisms.
Deformations of the natural semi-direct product
$S\# G$, the {\em skew group algebra},
include many compelling and influential algebras.
For example, when $G$ acts linearly on a finite dimensional,
complex vector space $V$, it induces an action on the symmetric algebra $S(V)$
(a polynomial ring).  Deformations of $S(V)\#G$ 
play a profound role in representation theory and connect diverse areas
of mathematics.
Graded Hecke algebras were originally
defined independently by Drinfeld~\cite{Drinfeld} 
and (in the special case of Weyl groups) by Lusztig~\cite{Lusztig3}.
These deformations of $S(V)\# G$ 
include symplectic reflection algebras (investigated by Etingof 
and Ginzburg~\cite{EtingofGinzburg} in the study of orbifolds) and
rational Cherednik algebras (introduced by Cherednik~\cite{Cherednik} 
to solve Macdonald's inner product conjectures).
Gordon~\cite{Gordon} used these algebras to prove a version of the
$n!$ conjecture for Weyl groups due to Haiman.

The deformation theory of an algebra is governed by its Hochschild cohomology
as a graded Lie algebra under the Gerstenhaber bracket.
Each deformation of the algebra arises from a (noncommutative) {\em Poisson structure}, 
that is, an element of Hochschild cohomology in degree 2 
whose Gerstenhaber square bracket is zero.  
Thus, a first step in understanding an algebra's deformation theory
is a depiction of the Gerstenhaber bracket.
In some situations, every noncommutative Poisson structure {\em integrates}, 
i.e., lifts to a deformation (e.g., see Kontsevich~\cite{Kontsevich}). 
Halbout and Tang~\cite{HalboutTang} investigate some of these structures
for algebras $C^{\infty}(M)\# G$ where $M$ is a real manifold, concentrating on the case
that $M$ has a $G$-invariant symplectic structure.

In this paper, we explore the rich algebraic structure of the Hochschild cohomology
of $S\# G$ with an eye toward deformation theory.
We describe the Gerstenhaber bracket on the Hochschild cohomology
of $S\# G$, for a general algebra $S$.
We then specialize to the case $S=S(V)$ to continue an analysis 
begun in two previous articles.
Hochschild cohomology is a Gerstenhaber algebra under two operations,
a cup product and a graded Lie (Gerstenhaber) bracket.  
In~\cite{paper2}, we examined the cohomology
$\HH^{\DOT}(S(V)\#G)$ as a graded algebra under cup product.
In this paper, we study its Gerstenhaber bracket and 
find noncommutative Poisson structures.
These structures catalog potential deformations of
$S(V)\#G$, most of which have not yet been explored.

For any algebra $A$ over a field $k$, 
both the cup product and Gerstenhaber bracket on 
the Hochschild cohomology $\HH^{\DOT}(A):=\Ext^{\DOT}_{A\ot A^{op}}(A,A)$
are defined initially on the bar resolution, a natural $A\ot A^{op}$-free
resolution of $A$. The cup product has another description as Yoneda
composition of extensions of modules, which can be transported to any
other projective resolution.  However, the Gerstenhaber bracket has 
resisted such a general description. 
In this paper, we use isomorphisms of cohomology 
which encode traffic between resolutions
to analyze $\HH^{\DOT}(S\# G)$ and unearth its Gerstenhaber bracket.

For the case $S=S(V)$, we use Demazure operators on one hand and
quantum partial differentiation on the other hand
to render the Gerstenhaber bracket and gain
theorems which predict its vanishing.
For example, the cohomology $\HHD(S\#G)$ breaks into a direct sum over $G$.
By invoking Demazure operators to implement automorphisms of cohomology,
we procure the following main result in Section~\ref{sec:degree-2}:

\medskip

\noindent {\bf Theorem~\ref{bracketalwayszero}.}
{\em The  bracket of any two elements $\alpha$, $\beta$ 
in $\HH^2(S(V)\#G)$ supported on group elements acting nontrivially
on $V$ is zero: $
[\alpha,\beta]=0\ .$}

\medskip

\noindent
As a consequence of the theorem,
the Gerstenhaber square bracket of every element $\HH^2(S(V)\#G)$ supported
off the kernel of the action of $G$ on $V$ defines a noncommutative 
Poisson structure.
In particular, if $G$ acts faithfully, any Hochschild $2$-cocycle
with zero contribution from the identity group element 
defines a noncommutative Poisson structure.

The cohomology $\HHD(S\#G)$ is graded not only by cohomological degree,
but also by polynomial degree.
In~\cite{SheplerWitherspoon}, we showed that
every {\em constant} Hochschild 2-cocycle (i.e., of polynomial degree 0)
defines a graded Hecke algebra.  Since graded Hecke algebras
are deformations of $S(V)\#G$ (see Section~\ref{previouspaper}), 
this immediately implies
that every constant Hochschild 2-cocycle has square bracket zero.
We articulate our automorphisms of cohomology using quantum
partial differentiation
to extend this result in Section~\ref{sec:zero}
to cocycles of {\em arbitrary} cohomological degree:

\medskip

\noindent{\bf Theorem~\ref{polydeg0}}
{\em The bracket of any two {\em constant} cocycles $\alpha$, $\beta$ 
in $\HHD(S(V)\#G)$ is zero: $[\alpha,\beta]=0\ .$}

\medskip

Supporting these two main ``zero bracket'' theorems, we give formulas for
calculating brackets at the cochain level in
Theorems~\ref{closedGB} and~\ref{general-circ}.
We apply these formulas to the abelian case in Example~\ref{abelian}
and express the bracket of 2-cocycles in terms of inner products
of characters.
We use this example in Theorem~\ref{nonzeroabeliansquarebracket} to show that
the hypotheses of Theorem~\ref{bracketalwayszero} can not be weakened and
that its converse is false for abelian groups.

We briefly outline our program.
In Section~\ref{Prelim}, we establish notation and definitions.
In Sections~\ref{Isomorphisms} and \ref{Liftingbracket}, we construct an 
explicit, instrumental
isomorphism between $\HHD(S\# G)$ and $\HHD(S,S\# G)^G$ and  
lift the Gerstenhaber bracket on $\HHD(S\#G)$ to an arbitrary resolution
used to compute $\HH^{\DOT}(S,S\# G)^G$.
In Sections~\ref{sec:koszul} and \ref{bracketsforS(V)}, 
we turn to the case $S=S(V)$ and examine isomorphisms of cohomology
developed in~\cite{paper1}.
We lift the Gerstenhaber bracket on $\HH^{\DOT}(S(V)\# G)$
from the bar resolution of $S(V)\#G$
to the Koszul resolution  of $S(V)$.
In Section~\ref{ebf}, we give our closed formulas for the bracket
in terms of quantum differentiation and recover the classical
Schouten-Nijenhuis bracket in case $G=\{1\}$.
These formulas allow us to characterize geometrically some cocycles with 
square bracket zero in Section~\ref{sec:zero}.  

Mindful of applications to 
deformation theory, we turn our attention to cohomological degree two
in Section~\ref{sec:degree-2}.
We prove the above Theorem~\ref{bracketalwayszero}
(which applies to cocycles of {\em arbitrary} polynomial degree)
and focus on abelian groups in Section~\ref{sec:abelian}.
We then apply our approach to the deformation theory of $S(V)\#G$
and graded Hecke algebras.
When a Hochschild 2-cocycle is a graded map of negative degree, corresponding
deformations are well understood and include graded Hecke algebras 
(which encompass symplectic reflection algebras and
rational Cherednik algebras).
In Section~\ref{previouspaper}, we explicitly present
graded Hecke algebras as deformations
defined using Hochschild cohomology, 
thus fitting our previous work on graded
Hecke algebras into a more general program to understand all
deformations of $S(V)\# G$.
We show how to render a constant Hochschild $2$-cocycle
into the defining relations of a graded Hecke algebra 
concretely and vice versa.
In fact, we use our results to explain
how to convert {\em any} Hochschild $2$-cocycle
into the explicit multiplication map 
of an infinitesimal deformation of $S(V)\# G$, allowing for exploration of a wide
class of algebras which include graded Hecke algebras as examples.

%%%%%%%%%%%%%%%%%%%%%%%%%%%%%%%%%%%%%%%%%%%%%%%%%%%%%%%%%%%%%%
%%%%%%%%%%%%%%%%%%%%%%%%%%%% SECTION %%%%%%%%%%%%%%%%%%%%%%%%%
%%%%%%%%%%%%%%%%%%%%%%%%%%%%%%%%%%%%%%%%%%%%%%%%%%%%%%%%%%%%%%
\section{Preliminary material}
\label{Prelim}

Let $G$ be a finite group. 
We work over the complex numbers $\CC$; all tensor products will be taken
over $\CC$ unless otherwise indicated.
Let $S$ be any $\CC$-algebra on which $G$ acts by automorphisms.
Denote by ${}^gs$ the result of applying an element $g$ of $G$ to
an element $s$ of $S$.
The {\bf skew group algebra} $S\# G$ is the vector space
$S\ot \CC G$ with multiplication given by
$$
  (r\ot g)(s\ot h)=r ( {}^{g}s) \ot gh
$$
for all $r,s$ in $S$ and $g,h$ in $G$.
We abbreviate $s\ot g$ by $s\overline{g}$ ($s\in S$, $g\in G$)
and $s\ot 1, \ 1\ot g$ simply by $s, \ \overline{g}$, respectively.
An element $g$ in $G$ acts on $S$ by an inner automorphism in $S\# G$: 
$ \ \overline{g} s (\overline{g})^{-1} = ( {}^gs)\overline{g} (\overline{g})
^{-1} = {}^gs$ for all $s$ in $S$.
We work with the induced group action on all maps
throughout this article:
For any map $\theta$ and element $h$ in $\Gl(V)$, 
we define $ {}^h\theta $ 
by $(^h\theta)(v) := \hphantom{\,}^h(\theta(^{h^{-1}}v))$ for all $v$.

%%%%%%%%%%%%%%%%%%%%%%%%%%%%%%%%%%%%%%%%%%%%%%%%%%%%%%%%%%%%%%%%%%%%%%%%%%%%%5
\subsection*{Hochschild cohomology and deformations}
The {\bf Hochschild cohomology} of a $\CC$-algebra $A$
with coefficients in an $A$-bimodule $M$
is the graded vector space
$\HH^{\DOT}(A,M)=\Ext^{\DOT}_{A^e}(A,M)$, where $A^e=A\ot A^{op}$ acts
on $A$ by left and right multiplication.
We abbreviate $\HH^{\DOT}(A)=\HH^{\DOT}(A,A)$.

One projective resolution of $A$ as an $A^e$-module is the {\bf bar
resolution}
\begin{equation}\label{barcomplex}
  \cdots\stackrel{\delta_3}{\longrightarrow} A^{\ot 4}
\stackrel{\delta_2}{\longrightarrow} A^{\ot 3}
\stackrel{\delta_1}{\longrightarrow} A^e
\stackrel{m}{\longrightarrow} A \rightarrow 0,
\end{equation}
where $\delta_p(a_0\ot\cdots\ot a_{p+1}) = \sum_{j=0}^p (-1)^j a_0
\ot\cdots\ot a_ja_{j+1}\ot\cdots\ot a_{p+1}$, and $\delta_0=m$ is
multiplication.
For each $p$, $\Hom_{A^e}(A^{\ot (p+2)},*)\cong \Hom_{\CC}(A^{\ot p}, *)$,
and we identify these two spaces of $p$-cochains
throughout this article.

The {\bf Gerstenhaber bracket} on Hochschild cohomology $\HH^{\DOT}(A)$
is defined at the chain level on the bar complex.
Let $f\in\Hom_{\CC}(A^{\ot p},A)$ and $f'\in\Hom_{\CC}(A^{\ot q},A)$.
The Gerstenhaber bracket $[f,f']$ in $\Hom_{\CC}(A^{\ot (p+q-1)}, A)$ 
is defined as
\begin{equation}\label{circ}
[f,f']=f\,\overline{\circ} f' - (-1)^{(p-1)(q-1)} f'\,\overline{\circ}f,
\end{equation}
where
$$
\begin{aligned}
  &f\overline{\circ}f'(a_1\ot\cdots\ot a_{p+q-1})\\
   &\hspace{.2cm}=
   \sum_{k=1}^p (-1)^{(q-1)(k-1)} f(a_1\ot\cdots\ot a_{k-1}\ot f'
  (a_k\ot\cdots\ot a_{k+q-1})\ot a_{k+q}\ot\cdots\ot a_{p+q-1})
\end{aligned}
$$
for all $a_1,\ldots, a_{p+q}$ in $A$.
This induces a bracket on 
$\HH^{\DOT}(A)$ under which it is a graded Lie algebra.
The bracket is compatible with the cup product on $\HH^{\DOT}(A)$, in the
sense that 
if $\alpha\in\HH^p(A)$, $\beta\in\HH^q(A)$, and $\gamma\in\HH^r(A)$, then
$$
[\alpha\smile\beta,\gamma] = [\alpha,\gamma]\smile\beta +
       (-1)^{p(r-1)} \alpha\smile [\beta,\gamma].
$$
Thus $\HH^{\DOT}(A)$ becomes a Gerstenhaber algebra.

Let $t$ be an indeterminate. 
A {\bf formal deformation} of $A$ is an associative $\CC[[t]]$-algebra
structure on formal power series $A[[t]]$ with multiplication determined
by 
$$
  a*b = ab + \mu_1(a\ot b)t + \mu_2(a\ot b)t^2+\cdots
$$ 
for all $a,b$ in $A$, where $ab=m(a\ot b)$ and the $\mu_i: A\ot A\rightarrow A$
are linear maps.
In case the above sum is finite for each $a,b$ in $A$, we may consider the
subalgebra  $A[t]$, a deformation of $A$ over the polynomial ring $\CC[t]$.
Graded Hecke algebras (see Section~\ref{previouspaper}) are examples
of deformations over $\CC[t]$.

Associativity of the product $*$ implies in particular that $\mu_1$ is a Hochschild
2-cocycle (i.e., $\delta_3^*(\mu_1)=0$) and that
the Gerstenhaber bracket $[\mu_1,\mu_1]$ is a coboundary 
(in fact, $[\mu_1,\mu_1]=2\delta_4^*(\mu_2)$, see~\cite[($4_2$)]{Gerstenhaber}).
If $f$ is a Hochschild 2-cocycle, then $f$ may or may not 
be the first multiplication map $\mu_1$ for some formal deformation.
The {\bf primary obstruction}
to lifting $f$ to a deformation is the Gerstenhaber bracket
$[f,f]$ considered as an element of $\HH^3(A)$.
A (noncommutative) {\bf Poisson structure} on $A$ is an element $\alpha$ in $\HH^2(A)$
such that $[\alpha,\alpha]=0$ as an element of $\HH^3(A)$.
The set of noncommutative Poisson structures includes the cohomology classes of the
first multiplication
maps $\mu_1$ of all deformations of $A$, and thus every deformation of $A$ defines a
noncommutative Poisson structure.
(For more details on Hochschild cohomology, see Weibel~\cite{Weibel}, and on 
deformations, see Gerstenhaber~\cite{Gerstenhaber};
Poisson structures for noncommutative algebras
are introduced in Block and Getzler~\cite{BlockGetzler} and Xu~\cite{Xu}.)
%%%%%%%%%%%%%%%%%%%%%%%%%%%%%%%%%%%%%%%%%%%%%%%%%555
%%%%%%%%%%%%%%%%%%%%%%%%%%%%%%%%%%%%%%%%%%%%%%%%%%%%%5
%%%%%%%%%%%%%%%%%%%%%%%%%%%%%%%%%%%%%%%%%%%%%%%%%%%%%%%%%5
\section{Hochschild cohomology of skew group algebras}
\label{Isomorphisms}

We briefly describe various isomorphisms used to compute 
Hochschild cohomology.
Let $S$ be any algebra upon which a finite group $G$ 
acts by automorphisms,
and let $A:=S\# G$ denote the resulting skew group algebra.
Let $\mathcal C$ be  a set of representatives of the conjugacy classes of $G$. 
For any $g$ in $G$, let $Z(g)$ be 
the centralizer of $g$.
A result of \c{S}tefan~\cite[Cor.\ 3.4]{Stefan} implies that
there is a $G$-action giving the first  
in a series of isomorphisms of graded vector spaces:
\begin{equation}\label{isos}
\begin{aligned}
  \HH^{\DOT}(S\# G) \ \ 
&\cong \ \ 
\HH^{\DOT}(S, S\#G)^G\\
\ \ &\cong \ \ 
 \Biggl(\,\bigoplus_{g\in G} \HH^{\DOT}(S, S\bar{g})
\Biggr)^{G}  \\
\ \ &\cong \ \ 
\bigoplus_{g\in{\mathcal C}}\HH^{\DOT}(S,S\bar{g})^{Z(g)}.
\end{aligned}
\end{equation}
The first isomorphism is in fact an isomorphism of graded algebras (under
the cup product); it follows from applying a spectral sequence.
We take a direct approach in this paper to gain results on the Gerstenhaber
bracket.
The second isomorphism results from decomposing 
the $S^e$-module $S\# G$ into the direct sum of components $S\,\overline{g}$.
The action of $G$ permutes these components via the conjugation action of $G$
on itself. The third isomorphism canonically projects onto a set of 
representative summands.

In the next theorem, we begin transporting
Gerstenhaber brackets on $\HHD(S\# G)$ to the other spaces
in~(\ref{isos})
by etching an {\em explicit isomorphism} between
$\HHD(S\# G)$ and $\HHD(S,S\# G)^G$.
First, some definitions.
We say that a projective resolution $P_{\DOT}$ of $S^e$-modules is
{\bf G-compatible} if $G$ acts on each term in the resolution
and the action commutes with the differentials.
Let $\rey$ be the {\bf Reynold's operator} on $\Hom_{\CC}(S^{\otimes \DOT},A)$,
$$\rey(\gamma) := \frac{1}{|G|}\sum_{g\in G}\, ^g\!\gamma.$$
We shall use the same notation $\rey$ to denote the analogous operator on
any vector space carrying an action of $G$.
Let $\rh^*$ be
the map ``move group elements far right, applying them along the way'':
Define
$\rh^*: \Hom_{\CC}(S^{\ot p},A)^G\rightarrow \Hom_{\CC}(A^{\ot p},A)$
by 
\begin{equation}\label{Theta-star}
  \rh^*(\kappa) (f_1\overline{g}_1\ot\cdots\ot f_p\overline{g}_p) =
  \kappa (f_1\ot {}^{g_1}f_2\ot \cdots\ot {}^{(g_1\cdots g_{p-1})}f_p) 
   \,\overline{g_1\cdots g_p}
\end{equation}
for all $f_1,\ldots,f_p$ in $S$ and $g_1,\ldots,g_p$ in $G$.
A similar map appears in~\cite{HalboutTang}.
In the next theorem, we give 
an elementary explanation for the appearance of $\rh^*$ in an explicit
isomorphism of cohomology; the proof shows that $\rh^*$ is induced
from a map $\rh$ at the chain level.

\vspace{3ex}
%%%%%%%%%%%%%%%%%%%%%%%%%%%%%%%%%%%%%%%%%%%%%%%%%%%%%%
\begin{remark}\label{invariantcohomology}
The notion of $G$-invariant cohomology here is well-formulated.
If $S$ and $B$ are algebras on which $G$ acts by automorphisms and $B$ is
an $S$-bimodule,
we may define $G$-invariant cohomology $\HH^{\DOT}(S,B)^G$ via any $G$-compatible
resolution.
This definition does not depend on choice of resolution.
Indeed, since $|G|$ is invertible in our field, any chain map between
resolutions can be averaged over the group to produce a $G$-invariant
chain map. By the Comparison Theorem~\cite[Theorem 2.2.6]{Weibel}, 
this yields not only an isomorphism
on cohomology $\HH^{\DOT}(S,B)$ arising from two different resolutions, but
a $G$-invariant isomorphism.
\end{remark}
%%%%%%%%%%%%%%%%%%%%%%%%%%%%%%%%%%%%%%%%%%%%%
\vspace{3ex}

We adapt ideas of Theorem~5.1 of~\cite{CaldararuGiaquintoWitherspoon} 
to our general situation. 
It was assumed implicitly in that theorem
that the resolution defining cohomology was $G$-compatible
and that the cochains fed into the given map were $G$-invariant.
We generalize that theorem while simultaneously
making it explicit.  
We apply a Reynold's operator to turn 
the cochains fed into the map $\rh^*$ invariant, as they might not
be so a priori.
This allows us to apply chain maps to arbitrary cocycles, in fact,
to cocycles that may not even represent invariant cohomology classes.

\vspace{2ex}
%%%%%%%%%%%%%%%%%%%%%%%%%%%%%%%%%%%%%%%%%%
\begin{thm}\label{S:CGW} 
Let $S$ be an arbitrary $\CC$-algebra upon which a finite group $G$
acts by automorphisms, and let $A=S\# G$.  The map 
$$
\begin{aligned}
\rh^*\circ \rey: \Hom_{\CC}(S^{\ot p}, A)
&\rightarrow
\Hom_{\CC}(A^{\ot p},A),\\
\end{aligned}
$$
given by 
$$
 \bigl((\rh^*\circ  \rey )(\beta)\bigr) (f_1\overline{g}_1\ot\cdots\ot f_p\overline{g}_p)=
 \frac{1}{|G|}\sum_{g\in G}\, ^g\! \beta
(f_1\ot {}^{g_1}\scooch f_2\ot\cdots\ot {}^{g_1\cdots g_{p-1}}f_p)\, \overline{g_1\cdots g_p}
$$
for all $f_1,\ldots, f_p\in S$,  $g_1,\ldots, g_p\in G$, and 
$\beta\in\Hom_{\CC}(S^{\ot p},A)$,
induces an isomorphism
$$ \HH^p(S,A)^G \stackrel{\sim}{\longrightarrow} \HH^p(A).  $$
\end{thm}
%%%%%%%%%%%%%%%%%%%%%%%%%%%%%%%%%%%%%%

We note that the following proof does not require the base field to be $\CC$,
only that $|G|$ be invertible in the field.

%%%%%%%%%%%%%%%%%%%%%%%%%%%%%%%%%%%%%%%%%%
\begin{proof}
We express every $G$-compatible resolution as
a resolution over a ring that succinctly absorbs the group and its action:
Let 
$$
\Delta :=\oplus_{g\in G}\ S\,\overline{g}\otimes S\,\overline{g^{-1}}\ ,
$$
a natural subalgebra of $A^e$ containing $S^e$.
Any $S^e$-module carrying an action of $G$ is naturally
also a $\Delta$-module:
Each group element $g$ acts as the ring element
$\overline{g}\ot \overline{g^{-1}}$ in $\Delta$. 
In fact, 
an $S^e$-resolution of a module is $G$-compatible if and only if
it is simultaneously a $\Delta$-resolution.

The bar resolution for $S$ is $G$-compatible, i.e.,
extends to a $\Delta$-resolution.
We define a different resolution of $S$ as a $\Delta$-module
that interpolates between the bar resolution of $S$ and the bar resolution of $A$.
For each $p\geq 0$, let
$$
   \Delta_p := \left\{ \sum f_0\overline{g}_0\ot\cdots\ot f_{p+1}\overline{g}_{p+1} \mid
   f_i\in S, \ g_i\in G \mbox{ and } g_0\cdots g_{p+1}=1\right\},
$$
a $\Delta$-submodule of $A^{\ot (p+2)}$.
Note that $\Delta_0=\Delta$, each $\Delta_p$ is a projective $\Delta$-module, 
and $\Delta_p$ is also 
a projective $S^e$-module by restriction. The complex
$$
  \cdots \stackrel{\delta_3}{\longrightarrow} \Delta_2
   \stackrel{\delta_2}{\longrightarrow} \Delta_1
   \stackrel{\delta_1}{\longrightarrow} \Delta_0 
   \stackrel{m}{\longrightarrow} S(V)\rightarrow 0
$$
(where $m$ is multiplication and $\delta_i$ is the restriction of the
differential on the bar resolution of $A$) is a $\Delta$-projective
resolution of $S$. 
By restriction it is also an $S^e$-projective resolution of $S$.
The inclusion $S^{\ot (p+2)}\hookrightarrow \Delta_p$ induces a restriction map
$$
 \Hom_{S^e}(\Delta_p, A) \rightarrow \Hom_{S^e}(S^{\ot (p+2)}, A)
$$
that in turn induces an automorphism on the cohomology $\HH^p(S,A)$.
The subspace of $G$-invariant elements of either of the above Hom spaces 
is precisely the subspace of $\Delta$-homomorphisms.
Therefore, a $G$-invariant element of $\Hom_{S^e}(S^{\ot (p+2)},A)$
identifies with an element of
$\Hom_{\Delta}(\Delta_p,A)$, and its application to an element of
$\Delta_p$ is found via the $\Delta$-chain map $\Delta_{\DOT}\rightarrow
S^{\ot (\DOT + 2)}$ given by
\begin{equation}\label{moveright2}
  f_0\overline{g}_0\ot\cdots\ot f_{p+1}\overline{g}_{p+1}\mapsto
    f_0\ot {}^{g_0}f_1\ot {}^{g_0g_1}f_2\ot \cdots \ot 
    {}^{(g_0g_1\cdots g_p)}f_{p+1}
\end{equation}
for all $f_0,\ldots,f_{p+1}\in S$ and $g_0,\ldots,g_{p+1}\in G$.
Finally, one may check that the following map is an
isomorphism of $A^e$-modules
$A^{\ot (p+2)}\stackrel{\sim}{\rightarrow} A^e\ot_{\Delta}\Delta_p$:
\begin{equation}\label{moveright1}
  f_0\overline{g}_0\ot\cdots\ot f_{p+1}\overline{g}_{p+1} \mapsto
  (1\ot \overline{g_0\cdots g_{p+1}})\ot (f_0\overline{g}_0\ot\cdots
   \ot f_{p+1} \overline{g}_{p+1}(\overline{g_0\cdots g_{p+1}})^{-1}).
\end{equation}
This gives rise to an isomorphism $\Hom_{\Delta}(\Delta_p,A)\cong
\Hom_{A^e}(A^{\ot (p+2)},A)$.
In fact, one may obtain the bar resolution of $A$
from the $\Delta$-resolution $\Delta_{\DOT}$ directly 
by applying the functor $A^e\ot_{\Delta}\ot - $;
the map (\ref{moveright1}) realizes the corresponding
Eckmann-Shapiro isomorphism on cohomology, $\Ext^{\DOT}_{\Delta}(S,A)
\cong \Ext^{\DOT}_{A^e}(A,A)$, at the chain level.
Let $\rh: A^{\ot (p+2)}\rightarrow A^e\ot_{\Delta}S^{\ot (p+2)}$ 
denote the composition of (\ref{moveright1}) and (\ref{moveright2}), that
is 
$$
   \rh(f_0\overline{g}_0\ot\cdots\ot f_{p+1}\overline{g}_{p+1}) 
   = (1\ot \overline{g_0\cdots g_{p+1}}) \ot f_0\ot {}^{g_0}f_1\ot {}^{g_0g_1}f_2\ot
   \cdots \ot {}^{(g_0g_1\cdots g_p)}f_{p+1}
$$
for all $f_0,\ldots,f_{p+1}\in S$ and $g_0,\ldots,g_{p+1}\in G$.
The induced map $\rh^*$ on cochains is indeed that defined by 
equation~(\ref{Theta-star}).
The above arguments show that this induced map $\rh^*$ gives an explicit isomorphism, at
the chain level, from $\HH^{p}(S,A)^G$ to $\HH^p(A,A)$.
\end{proof}
\vspace{1ex}

%%%%%%%%%%%%%%%%%%%%%%%%%%%%%%%%%%%%%%%%%%%%%%%%%%%%%%%%%
%%%%%%%%%%%%%%%%%%%%%%%%%%%%%%%%%%%%%%%%%%%%%%%%%%%%%%%5
%%%%%%%%%%%%%%%%%%%%%%%%%%%%%%%%%%%%%%%%%%%%%%%%%%%%%%%5
\section{Lifting brackets to other resolutions}\label{Liftingbracket}

The Gerstenhaber bracket on the Hochschild cohomology of an algebra 
is defined using the bar resolution of that algebra.
But cohomology is often computed using some other projective
resolution.  One seeks to express the induced Gerstenhaber bracket
on any other projective resolution giving cohomology.  
Again, let $S$ be any algebra upon which a finite group $G$ 
acts by automorphisms,
and let $A:=S\# G$ denote the resulting skew group algebra.
In this section, we lift the Gerstenhaber bracket on $\HHD(A)$ 
to other resolutions used to compute cohomology. 
This task is complicated by the fact that we compute the cohomology
$\HHD(A)$ as the space $\HHD(S,A)^G$ (using the isomorphisms of~(\ref{isos})).  Hence,
we consider alternate resolutions of $S$, not $A$.
 
Suppose the Hochschild cohomology of $S$ has been determined using
some $S^e$-projective resolution $P_{\DOT}$ of $S$.
We assume that $P_{\DOT}$ is $G$-compatible, so that
$P_{\DOT}$ defines the $G$-invariant cohomology $\HHD(S,A)^G$.
Let $\Psi$ and $\Phi$ be chain maps
between the bar resolution and $P_{\DOT}$, i.e.,
\begin{eqnarray*}
   \Psi_p & : & S^{\ot (p+2)}\rightarrow P_p\, ,\\
   \Phi_p &:& P_p \rightarrow S^{\ot (p+2)}\, ,
\end{eqnarray*}
for all
$p\geq 1$, and the following diagram commutes:
\begin{equation}
\label{diagram}
\xymatrix{
\cdots \ar[r] & S^{\ot {(p+2)}}\ar[r]^{\delta_{p}}\ar@<-2pt>[d]_{\Psi_{p}} 
               & S^{\ot {(p+1)}}\ar[r]\ar@<-2pt>[d]_{\Psi_{p-1}} 
               & \cdots                  \\
\cdots \ar[r] & P_{p} \ar[r]^{d_{p}}\ar@<-2pt>[u]_{\Phi_{p}}
               & P_{p-1} \ar[r]\ar@<-2pt>[u]_{\Phi_{p-1}}
               & \cdots                 \ .
}
\end{equation}
Let $\Psi^*$ and $\Phi^*$ denote the maps induced by application 
of the functor $\Hom_{S^e}(-,A)$.

%%%%%%%%%%%%%%%%%%%%%%%5
\vspace{2ex}
\begin{lemma}\label{invariantisos}
The cochain maps $\Psi^*$ and $\Phi^*$ induce $G$-invariant, inverse
isomorphisms on the cohomology $\HHD(S,A)$.
\end{lemma}
\vspace{1ex}
%%%%%%%%%%%%%%%%%%%%%%%%
\begin{proof}
Since both $P_{\DOT}$ and the bar resolution are $G$-compatible, ${}^g\Psi^*$ and ${}^g\Phi^*$
are chain maps for any $g$ in $G$.
Such chain maps are unique up to chain homotopy equivalence,
and so induce the same maps on cohomology.
\end{proof}
\vspace{1ex}

We represent an element of $\HH^{\DOT}(S,A)$ at the chain level by
a function $f:P_p\rightarrow A$.
By Theorem~\ref{S:CGW},
the corresponding function from $A^{\ot (p+2)}$ to $A$ is given by
$\rh^*\circ\rey \circ\Psi_p^*(f)$.
Lemma~\ref{invariantisos} thus implies:

%%%%%%%%%%%%%%%%%%%%%%%%%%%%%%%%%%%%%%%%%%%%%%%%%%%%%%%%
\vspace{2ex}
\begin{thm}\label{general-S}
Let $A=S\# G$.
Let $P_{\DOT}$ be any $G$-compatible $S^e$-resolution 
of $S$ and
let $\Psi$ be a chain map from the bar resolution of $S$ to $P_{\DOT}$.
The map $$\rh^*\circ \rey\circ\Psi^* : 
\Hom_{S^e}(P_p,A)\rightarrow \Hom_{\CC}(A^{\ot p},A)$$
induces an isomorphism $\HH^p(S,A)^G\stackrel{\sim}{\longrightarrow} \HH^p(A)$.
\end{thm}
%%%%%%%%%%%%%%%%%%%%%%%%%%%%%%%%%%%%%%%%%%%%%%%%%%%%%%%%%%%%%%%%%%%5

We now describe the inverse isomorphism.
The inclusion map $S \hookrightarrow A$ induces a restriction map
$$
  \resAS : \Hom_{\CC}(A^{\ot \DOT}, A)\rightarrow \Hom_{\CC}(S^{\ot\DOT},A).
$$

%%%%%%%%%%%%%%%%%%%%%%%%%%%%%%%%%%%%%5
\vspace{2ex}
\begin{thm}\label{general-S2}
Let $A=S\# G$. 
Let $\Psi, \Phi$ be any chain maps converting between the bar resolution of $S$
and $P_{\DOT}$, as above.
The maps 
$$
  \Phi^*\circ\resAS \ \ \mbox{ and } \ \ 
   \rh^*\circ\rey\circ\Psi^*\ 
$$
$$
\begin{xy}
\xymatrix{ 
\text{on cochains, }
 \Hom_{\CC}(A^{\ot p},A) \  
\ar@<2pt>[rr]^{\hspace{0cm}\Phi^*\,\circ\,\res\hspace{1cm}} &&
\ \Hom_{S^e}(P_{p},A)\, ,
 \ \ar@<2pt>[ll]^{\hspace{0cm}\Theta^*\,\circ\,\rey\,\circ\,\Psi^*\hspace{.9cm}} 
\text{induce inverse isomor-}
}
\end{xy}
$$
phisms on cohomology,
$$\HHD(A)\stackrel{\sim}{\longleftrightarrow}
\HHD(S,A)^G\ .$$
\end{thm} 
%%%%%%%%%%%%%%%%%%%%%%%%%%%%%
\vspace{1ex}
\begin{proof}
By Theorem~\ref{general-S}, $\Theta^*\circ\rey\circ\Psi^*$ induces an isomorphism from 
$\HH^p(S,A)^G$ to $\HH^p(A)$.
We show that $\Phi^*\circ \resAS$ defines an inverse map on
cohomology.

By functoriality of Hochschild cohomology (see~\cite[\S1.5.1]{Loday}),
$\resAS$ induces a linear map:
$$
\resAS: \HH^{p}(A) \rightarrow \HH^{p}(S,A).
$$
We saw in the proof of Theorem~\ref{S:CGW} that
the  Eckmann-Shapiro isomorphism realized at the chain
level in (\ref{moveright1}) 
induces a map from $\Hom_{A^{e}}(A^{\ot(p+2)},A)$
to $\Hom_{S^e}(\Delta_p,A)^G$.
The restriction of (\ref{moveright1}) to $S^{\ot(p+2)}$
is essentially the identity map, and it induces 
a map to $\Hom_{S^e}(S^{(p+2)}, A)^G$
agreeing with $\resAS$.
Hence the image of $\resAS$ is $G$-invariant.
Note that $\resAS\circ \rh^* = 1$ on cochains, and thus on cohomology.
By the proof of Theorem~\ref{S:CGW}, $\rh^*$ is invertible on cohomology, 
and so we conclude that $\resAS=(\rh^*)^{-1}$
and $\resAS$  also yields an isomorphism on cohomology:
$$
 \resAS:  \HH^{\DOT}(A)\stackrel{\sim}{\longrightarrow}
 \HH^{\DOT}(S,A)^G.
$$

Lemma~\ref{invariantisos} then implies
that $\Phi^*\circ\res$ induces a well-defined map on cohomology
with $G$-invariant image.
The fact that it is  inverse to $\Theta^*\circ\rey\circ\Psi^*$ follows from the observations that
$\rey$ commutes with $\Phi^*$ (as $\Phi^*$
is invariant), $\Phi^*$ and $\Psi^*$ are 
inverses on cohomology, and
$\rey$ is the identity on $G$-invariants.
\end{proof}
\vspace{1ex}

%%%%%%%%%%%%%%%%%%%%%%%%%
\vspace{3ex}
\begin{remark}\label{replace}
We may replace $\Psi^*$ in Theorems~\ref{general-S} and~\ref{general-S2}
above by any other cochain map (with the same domain and range),
provided that map induces an automorphism inverse to $\Phi^*$ on the 
cohomology $\HHD(S,A)$, even if it is not induced by a map on the original
projective resolution.
(As $\Phi^*$ is $G$-invariant on cohomology, so too is its inverse.)
\end{remark}
\vspace{3ex}
%%%%%%%%%%%%%%%%%%%%%%%%%%%%%%%%%%%%%%%%%%%%%%%%%%
\begin{remark}\label{RR}
Some comments are in order before we 
give a formula for the Gerstenhaber bracket.
By Theorem~\ref{general-S2}, a cocycle in $\Hom_{\CC}(A^{\ot p},A)$
is determined by its values on
$S^{\otimes p}$, and so we may compute the bracket of two cocycles by
determining how that bracket acts on elements in $S^{\ot p}$.
The ring $S$ embeds canonically in $A$, and we identify this ring with its
image when convenient.
We also note that, for the purpose of using Theorem~\ref{general-S2},
it suffices to start with a (not necessarily $G$-invariant) cocycle $\alpha$ in
$\Hom_{S^e}(P_{\DOT},A)$ and apply $\rh^*\circ\rey\circ\Psi^*$, 
since $\rey\circ\Psi^*(\alpha)$
is cohomologous to $\rey\circ\Psi^*\circ\rey(\alpha)$:
Indeed, observe that
$\rey(\Psi^*(\rey(\alpha)))=\rey(\rey(\Psi^*)(\alpha))$ 
and that
$ \rey({\rey}(\Psi^*)(\alpha)) \sim \rey(\Psi^*(\alpha))$
since $\Psi^*$ and ${\rey}(\Psi^*)$ are both chain maps.
\end{remark}
\vspace{3ex}
%%%%%%%%%%%%%%%%%%%%%%%%%%%%%%%%%%%%%%%%%%%%%%%%%%

We now lift the Gerstenhaber
bracket on $\HHD(A)$ to the cohomology $\HHD(S, A)^G$
expressed in terms of any resolution of $S$.
Again, we assume the resolution is $G$-compatible, or else
it may not define $\HHD(S,A)^G$. 
Theorem~\ref{general-S2} implies the following formula for the graded
Lie bracket at the cochain level, on cocycles in $\Hom_{S^e}(P_{\DOT},A)^G$.
We use Definition~\ref{circ} of the bracket on $\HHD(A)$ at the cochain level.
%%%%%%%%%%%%%%%%%%%%%%%%%%%%%%%%%%%%%%%%%%%%%%%%%%%%%%%%%%%
\vspace{2ex}
\begin{thm}\label{general-S3}
Let $P_{\DOT}$ be any $G$-compatible $S^e$-resolution 
of $S$.  
The isomorphism $\HH^{\DOT}(A) \cong \HHD(S,A)^G$ induces
the following bracket on $\HHD(S,A)^G$
(expressed via $P_{\DOT}$):
For $\alpha, \beta$ in $\Hom_{S^e}(P_{\DOT},A)$
representing $G$-invariant cohomology classes,
the cohomology class of $[\alpha,\beta]$ is represented
at the cochain level by
$$
  \frac{1}{|G|^2} \sum_{a,b\in G} \Phi^* [ \rh^* \,{}^a(\Psi^*\alpha),
  \rh^* \,{}^b(\Psi^*\beta)],
$$
where the bracket $[\ \ ,\ \,]$ on the right side is the
Gerstenhaber bracket on $\HH^{\DOT}(A)$.
\end{thm} 

\begin{proof}
We apply the inverse isomorphisms of Theorem~\ref{general-S2} to lift the
bracket from the bar complex of $A$ to $P_{\DOT}$:
The above formula gives the resulting bracket 
$ \Phi^* \res\, [ \Theta^*\rey\Psi^* (\alpha), \Theta^*\rey\Psi^*(\beta) ]$.
Note that the restriction map $\resAS$ is not needed in the formula
since the output of $\Phi$ is automatically in $S^{\ot \DOT}$.
Also note
that it is not necessary to apply the Reynolds operator to $\alpha$ and
$\beta$ before applying $\Psi^*$ since we are interested only in the bracket
at the level of cohomology (see Remark~\ref{RR}).
\end{proof}
\vspace{1ex}

%%%%%%%%%%%%%%%%%%%%%%%%%%%%%%%%%%%%%%%%%%%%%%%%%%%%%%%%%%%%%%
%%%%%%%%%%%%%%%%%%%%%%%%%%%%%%%%%%%%%%%%%%%%%%%%%%%%%%%%%%%%%%
\section{Koszul resolution}\label{sec:koszul}

Let $G$ be a finite group and $V$ a 
(not necessarily faithful) $\CC G$-module 
of finite dimension $n$.
The Hochschild cohomology of the skew group algebra 
$\HH^{\DOT}(S(V)\# G)$ is computed using the Koszul resolution 
of the polynomial ring $S(V)$ while the cup product and 
Gerstenhaber bracket are defined on the bar resolution
of $S(V)\# G$.  
In this section, we use machinery developed in Section~\ref{Liftingbracket}
to translate between spaces and between resolutions.

First, some preliminaries.
We denote the image of $v$ in $V$ under the action of 
$g$ in $G$ by ${}^gv$.  
Let $V^*$ denote the contragredient (or dual) representation.
For any basis $v_1,\ldots,v_n$ of $V$, let
$v_1^*,\ldots,v_n^*$ be the dual basis of $V^*$.
Denote the set of $G$-invariants in $V$ by
$V^G:=\{v \in V: {}^{g}v=v \text{ for all } g\in G\}$ 
and the $g$-invariant subspace of $V$ by
$V^g := \{v\in V: {}^gv=v\}$ for any $g$ in $G$. 
Since $G$ is finite, we may assume $G$ acts by isometries on $V$ (i.e.,
$G$ preserves a Hermitian form $\langle \ , \ \rangle$).
If $h$ lies in the centralizer $Z(g)$, then $h$ preserves both $V^g$ and its orthogonal 
complement $(V^g)^\perp$ (defined with respect to the Hermitian form).
We shall frequently use the observation that $(V^g)^{\perp} = \Ima (1-g)$.

%%%%%%%%%%%%%%%%%%%%%%%%%%%%%%%%%%%%%%%%%%%%%%%%%%%%%%%%%%%%%%
The {\bf Koszul resolution} $K_{\DOT}(S(V))$ is given by 
$K_0(S(V))=S(V)^e$ and
\begin{equation}\label{koszul-res2}
  K_{p}(S(V)) = S(V)^e \ot \Wedge^{p}(V)
\end{equation}
for $p\geq 1$,
with differentials defined by
\begin{equation}\label{koszul-diff}
d_p(1\ot 1\ot v_{j_1}\wedge\cdots\wedge v_{j_p}) = 
   \sum_{i=1}^p (-1)^{i+1} (v_{j_i}\ot 1 - 1\ot v_{j_i})\ot
   (v_{j_1}\wedge\cdots\wedge \hat{v}_{j_i}\wedge\cdots\wedge v_{j_p})
\end{equation}
for any choice of basis $v_1,\ldots, v_n$ of $V$
(e.g., see Weibel~\cite[\S4.5]{Weibel}).
We identify $\Hom_{\CC}(\Wedge^pV,S(V)\overline{g})$
with $S(V)\overline{g}\ot\Wedge^p V^*$ for each $g$ in $G$.

We fix the set of {\bf cochains} arising from the Koszul resolution
(from which the cohomology classes emerge)
as vector forms on $V$ 
tagged by group elements:
Let 
\begin{equation}\label{C-cochains}
C^{\DOT}=\bigoplus_{g\in G} C_g^{\DOT},
\quad\text{where}\quad 
C_g^p :=  S(V)\bar{g} \otimes \Wedge^{p} V^*
\quad\text{ for }\ g\in G.
\end{equation}
We refer to $C_g^{\DOT}$ as the set of cochains 
{\bf supported on} $g$.
Similarly, if $X$ is a subset of $G$, we
set $C^{\DOT}_X := \oplus_{g\in X}C^{\DOT}_g$, 
the set of cochains
{\bf supported on} $X$. 
Note that group elements permute the summands 
of $C^{\DOT}$ via the conjugation action of $G$ on itself.

We apply the notation of Section~\ref{Liftingbracket}
to the case $S=S(V)$.
Let $P_{\DOT}$ be the Koszul resolution (\ref{koszul-res2})
of $S(V)$.
Since the bar and Koszul complexes of $S(V)$ are both $S(V)^e$-resolutions,
there exist chain maps $\Phi$ and $\Psi$ between the two complexes,
\begin{eqnarray*}
   \Psi_p & : & S(V)^{\ot (p+2)}\rightarrow S(V)^e\ot \Wedge^p V,\\
   \Phi_p &:& S(V)^e\ot \Wedge^p V\rightarrow S(V)^{\ot (p+2)},
\end{eqnarray*}
for all $p\geq 1$, such that Diagram~\ref{diagram} commutes.
Let $\Phi$ be the canonical inclusion of the Koszul resolution~(\ref{koszul-res2})
into the bar resolution~(\ref{barcomplex}):
\begin{equation}\label{phik}
  \Phi_p(1\ot 1\ot v_{j_1}\wedge\cdots\wedge v_{j_p})
   = \sum_{\pi\in \Sym_p}\sgn(\pi)\ot v_{j_{\pi(1)}}\ot\cdots\ot
    v_{j_{\pi(p)}}\ot 1
\end{equation}
for all $v_{j_1},\ldots,v_{j_p}$ in $V$.
(See~\cite{paper1} for an explicit chain map $\Psi$ in this case.)
We obtain the following commutative diagram
of induced cochain maps:
\begin{equation}\label{big-diagram}
\xymatrix{
 \Hom_{\CC}(S(V)^{\ot p}, A) \ar[r]^{\delta^*}\ar@<2pt>[d]^{\Phi^*_p}
     & \Hom_{\CC}(S(V)^{\ot (p+1)}, A)\ar@<2pt>[d]^{\Phi^*_{p+1}}\\
   C^p\ar[r]^{d^*}  \ar@<2pt>[u]^{\Psi^*_p}
    & C^{p+1}\ar@<2pt>[u]^{\Psi^*_{p+1}}\ \ \ .
}
\end{equation}
Note that both the Koszul and bar resolutions are $G$-compatible.
Using Lemma~\ref{invariantisos}, we identify $\Phi^*$
and $\Psi^*$ with their restrictions to 
$$\HHD(S(V),S(V)\#G)^G\ .$$

Given any basis $v_1, \ldots, v_n$ of $V$, let $\partial /\partial v_i$ denote the usual
partial differential operator with respect to $v_i$.
In addition, given a complex number $\epsilon\neq 1$, we 
define the $\epsilon$-{\bf quantum partial
differential operator} with respect to
$v:=v_i$ as the scaled Demazure
(BGG) operator 
$\del_{v, \epsilon}:S(V)\rightarrow S(V)$
given by
\begin{equation}\label{qpd}
\del_{v, \epsilon}(f)\ = \ (1-\epsilon)^{-1}\ \frac{f-\, ^s\hspace{-.5ex}f}{v}
\ = \ \frac{f-\, ^s\hspace{-.5ex}f}{v-\, ^sv}\ ,
\end{equation}
where $s\in\Gl(V)$ is the reflection whose matrix with respect to the 
basis $v_1, \ldots, v_n$ is 
$\diag(1,\ldots, 1, \epsilon, 1, \ldots, 1)$
with $\epsilon$ in the $i$th slot.
 Set $\del_{v,\ep}= \del / \del v$
when $\ep=1$.
The operator $\del_{v,\epsilon}$ 
coincides with the usual definition of
quantum partial differentiation:
$$
\del_{v_1,\epsilon}(v_1^{k_1}v_2^{k_2}\cdots v_n^{k_n})
=[k_1]_{\epsilon} \  v_1^{k_1-1}v_2^{k_2}\cdots v_n^{k_n}\, ,
$$
where $[k]_{\epsilon}$ is the quantum integer
$
   [k]_{\epsilon} := 1 + \epsilon + \epsilon^2 +\cdots + \epsilon^{k-1}.
$

Next, we recall an explicit map $\ta$, 
involving quantum differential operators, that 
replaces $\Psi^*$:
\begin{equation}\label{big-diagram2}
\xymatrix{
 \Hom_{\CC}(S(V)^{\ot p}, A) \ar[r]^{\delta^*}\ar@<2pt>[d]^{\Phi^*_p}
     & \Hom_{\CC}(S(V)^{\ot (p+1)}, A)\ar@<2pt>[d]^{\Phi^*_{p+1}}\\
  C^p\ar[r]^{d^*}  \ar@<2pt>[u]^{\ta_p}
    & C^{p+1}\ar@<2pt>[u]^{\ta_{p+1}}\ \ \ .
}
\end{equation}
For each $g$ in $G$, fix a basis $B_g=\{v_1,\ldots, v_n\}$
of $V$ consisting of eigenvectors of $g$
with corresponding eigenvalues $\epsilon_1, \ldots, \epsilon_n$.
Decompose $g$ into reflections according to this basis: 
Let $g=s_1\cdots s_n$ where each $s_i$ in $\GL(V)$ is the reflection 
(or the identity) defined by
$s_i(v_j) = v_j$ for $j\neq i$ and $s_i(v_i)=\epsilon_i v_i$.
Let $\del_{i}:=\del_{v_i, \epsilon_i}$, 
the quantum partial derivative 
with respect to $B_g$.

%%%%%%%%%%%%%%%%%%%%%%%%%%%%%%%%%%%%%%%%%%%%%%%%%%%%%%%
\vspace{3ex}
\begin{defn}
\label{tau}
We define a linear map $\ta$ from 
the dual Koszul complex to the dual bar complex
for $S(V)$ with coefficients in $A:=S(V)\# G$,
$$
\ta_p:\ \ 
C^p 
\rightarrow
\Hom_{\CC}(S(V)^{\ot p}, A) \ .
$$
Fix $g$ in $G$ with basis $B_g$ of $V$ as above.
Let
$
\alpha = f_g \overline{g}\otimes v_{j_1}^*\wedge\cdots\wedge v_{j_p}^*$
lie in $C^p$
with 
$f_g$ in $ S(V)$ and $1\leq j_1<\ldots<j_p \leq n$.
Define $\ta(\alpha):S(V)^{\ot p}\rightarrow S(V)\bar{g}$ 
by
$$
\begin{aligned}
\ta(\alpha)(& f_1\otimes \cdots\otimes f_p )
&= \Biggl(\
\prod_{k=1,\ldots, p}\
^{s_{1}s_{2}\cdots s_{j_{k}-1}} 
(\del_{j_{k}}f_k) \Biggr) f_g\overline{g}
\ .
\end{aligned}
$$
Then $\ta$ is a cochain map (see~\cite{paper1}) and thus
induces an endomorphism $\ta$ of cohomology 
$\HH^{\DOT}(S(V), A)$.
Denote the restriction 
of $\ta$ to the $g$-component of $C^{\DOT}$ and
of $\HH^{\DOT}(S(V),A)\cong \oplus_{g\in G}\, \HH^{\DOT}(S(V),S(V)\bar{g})$
by $\ta_g = \ta_{g, B_g}$ so that
$$\ta = \bigoplus_{g\in{G}} \ta_g.
$$  
\end{defn}
%%%%%%%%%%%%%%%%%%%%%%%%%%%%%%%%%%%%%%
\vspace{3ex}
In formulas, we wish to highlight group elements explicitly instead
of leaving them hidden in the definition of $\ta$, and
thus we find it convenient to define a version of $\ta$ untagged by group elements:
\vspace{3ex}
\begin{defn}
\label{ttau}
Let $\proj_{G\rightarrow 1}: S(V)\# G\rightarrow S(V)$ 
be the projection map that drops group element 
tags: $f\overline{g}\mapsto f$ for all  $g$ in $G$ and $f$ in $S(V)$.  Define
$$
\ttau\ :=\ \proj_{G\rightarrow 1} \circ \ta:\ C^{\DOT} \rightarrow \bigoplus_{g\in G} 
\Hom_{\CC}(S(V)^{\ot \DOT},S(V))\ .
$$
\end{defn}
%%%%%%%%%%%%%%%%%%%%%%%%%%%%%%%%%%%%%%%%%%%%%%%%%%%%%%%%%%%%
\vspace{3ex}

We shall use the following observation often.

%%%%%%%%%%%%%%%%%%%%%%%%%%%%%%%%%%%%%%%%%%%%%%%%%%5
\vspace{3ex}
\begin{remark}
\label{zeroterm}
For the fixed basis $B_g=\{v_1,\ldots, v_n\}$ 
and $\alpha=f_g \overline{g}\otimes v_{j_1}^*\wedge\cdots\wedge v_{j_p}^*$ 
in $C^p_g$ (with $j_1<\ldots< j_p$), note that 
$$
\ta(\alpha)(v_{i_1} \otimes\cdots\ot v_{i_p})=0
\quad\quad\quad
\text{unless } i_1=j_1, \ldots, i_p=j_p\ .
$$
Generally, 
$
\ta(\alpha)(f_1\ot\cdots\ot f_p)=0
$
whenever ${\frac{\del}{\del v_{j_k}}(f_k)} =0$ for some $k$. 
\end{remark}
\vspace{3ex}

The following proposition from~\cite{paper1} provides a cornerstone for our computations.

\vspace{2ex}
%%%%%%%%%%%%%%%%%%%%%%%%%%%%%%%%%%%%%%%%%%%%%%%%%%%%%%5
\begin{prop}\label{phi-inverse} 
The map $\ta$ induces an isomorphism
on the Hochschild cohomology of the skew group algebra
$$\HH^{\DOT}(S(V)\# G)
\cong
\left(\,\bigoplus_{g\in G}\, \cohg  \right)^G\ .$$
Specifically, $\ta$ and $\Phi^*$ are
$G$-invariant, inverse isomorphisms on cohomology
converting between expressions in terms of the Koszul
resolution and the bar resolution.
\end{prop}
\vspace{2ex}
%%%%%%%%%%%%%%%%%%%%%%%%%%%%%%%%%%%%%%%%%%%%%%%%%%%%%%%%%%%%%%%%%%%%
In fact,
the map $\ta$ is easily seen to be a right inverse 
to $\Phi^*$ on cochains, not just on cohomology; see~\cite[Prop.\ 5.4]{paper2}.

%%%%%%%%%%%%%%%%%%%%%%%%%%%
\vspace{3ex}
\begin{remark}\label{dependsonbasis}
The {\em cochain} map $\ta=\bigoplus_{g\in G}\ta_{g,B_g}$ 
depends on the choices of
bases $B_g$ of eigenvectors of $g$ in $G$,
but $\ta$ induces an automorphism on cohomology 
which does {\em not} depend on the choice of basis. 
(This follows from Proposition~\ref{phi-inverse} above, as $\Phi$ is independent
of basis.)
\end{remark}
\vspace{3ex}
%%%%%%%%%%%%%%%%%%%%%%%%%%%%%%%%%%%%%%%%%%%%%%%%%%%%%%%%%%%%%%%%

%%%%%%%%%%%%%%%%%%%%%%%%%%%%%%%%%%%%%%%%%%%%%%%%%%%%%%%%%%%%%%%%%
%%%%%%%%%%%%%%%%%%%%%%%%%%%%%%%%%%%%%%%%%%%%%%%%%%%%%%%%%%%%%%%%%%
%%%%%%%%%%%%%%%%%%%%%%%%%%%%%%%%%%%%%%%%%%%%%%%%%%%%%%%%%%%%%%%%%%

\section{Brackets for polynomial skew group algebras}
\label{bracketsforS(V)}

Let $G$ be a finite group and let $V$ be a finite dimensional $\CC G$-module. 
We apply the previous results to determine the Gerstenhaber bracket
of $\HHD(S(V)\# G)$ explicitly.
We lift the Gerstenhaber bracket
on $\HH^{\DOT}(S(V)\# G)$, which is defined via the bar complex
of $S(V)\# G$, to the cohomology $\HHD(S(V), S(V)\#G)^G$,
which is computed via the Koszul complex of $S(V)$.

The next theorem allows one to replace the Hochschild cohomology of $S(V)\#G$ 
with a convenient space of forms. 
Set
\begin{equation}\label{H,Hg}
\begin{aligned}
    H^{\DOT}_g \ \    
     & \ \ :=\ S(V^g)\overline{g}\ot  \Wedge^{\DOT - \codim V^g}(V^g)^*
     \ot \Wedge^{\codim V^g}((V^g)^{\perp})^*\quad\text{and} \\
   H^{\DOT} \ \ & \ \ :=\  \bigoplus_{g\in G}H_g^{\DOT}\ ,
\end{aligned}
\end{equation}
where a negative exterior power is defined to be 0.
We regard these spaces as subsets of $C^{\DOT}$, the space of cochains 
arising from the Koszul complex (see (\ref{C-cochains})),
after making canonical identifications.
Note that for each $g$ in $G$, 
the centralizer $Z(g)$ acts diagonally on the tensor factors
in $H_g^{\DOT}$. 
The third factor $\Wedge^{\codim V^g}((V^g)^{\perp})^*$ of $H_g^{\DOT}$
is a vector
space of dimension one with a possibly nontrivial $Z(g)$-action.
The theorem below is due to
Ginzburg and Kaledin~\cite[(6.4)]{GinzburgKaledin} for 
faithful group actions in a more general geometric setting. 
See also Farinati~\cite[Section 3.2]{Farinati} in the algebraic setting.
Our formulation is
from~\cite[(3.3) and (3.4)]{SheplerWitherspoon}.

\vspace{2ex}
%%%%%%%%%%%%%%%%%%%%%%%%%%%%%%%%%%%%%%%5
\begin{thm} 
\label{Decomposition}
Set $A:=S(V)\# G$.
Cohomology classes arise as tagged vector forms:
\begin{itemize}
\item
$H^{\DOT}$ is a set of cohomology class representatives
for $\HHD(S(V), A)$ arising from the Koszul resolution:
The inclusion map $H^{\DOT}\hookrightarrow C^{\DOT}$ induces
an isomorphism $H^{\DOT} \cong \HHD(S(V), A)$.
\item 
$\HH^{\DOT}(A) \ \cong \ 
\HHD(S(V), A)^G \ \cong\
 (H^{\DOT})^G\ . $
\end{itemize}
\end{thm}
%%%%%%%%%%%%%%%%%%%%%%%%%%%%%%%%%%%%%%%%%%%%%%%%

\vspace{3ex}
\begin{remark}\label{degree-two}
The only contribution to degree two
cohomology $\HH^2(S(V)\#G)$ comes from group elements acting either trivially
or as bireflections:
$$H^2_g\neq 0
\quad\text{implies}\quad
V^g=V \ \ \text{or} \ \ \codim V^g=2\ .
$$
Indeed, by Definition~\ref{H,Hg},
$(H^{\DOT})^G\cong \oplus_{g\in {\mathcal C}}(H^{\DOT}_g)^{Z(g)}$
where $\mathcal C$ is a set of representatives of the conjugacy classes of $G$.
Since each group element $g$ in $G$ centralizes itself, 
the determinant of $g$ on $V$ is necessarily $1$ 
whenever $(H^{\DOT}_g)^{Z(g)}$ is nonzero
(since $g$ acts on 
$\Wedge^{\codim V^g}((V^g)^{\perp})^*$ by the inverse of its determinant).
In fact, if $g$ acts nontrivially on $V$, it does not contribute
to cohomology in degrees 0 and 1.
\end{remark}
\vspace{3ex}

Our results from previous sections allow us to realize the isomorphism
of Theorem~\ref{Decomposition} explicitly at the chain level in the
next theorem, which will be used extensively in our bracket
calculations: 

%%%%%%%%%%%%%%%%%%%%%%%%%%%%%%%%%%%%%%%%%%%%%%%%%%%%%%%%%
\vspace{2ex}
\begin{thm}\label{S:CGW2}
Let $A=S(V)\# G$. The map
$$
\begin{aligned}
\gamm : =\rh^*\circ \rey \circ \ta:\ \ \ 
   C^{p}
&\ \ {\rightarrow} \ \ 
\Hom_{\CC}(A^{\ot p},A)\\
 \alpha 
&\ \ \mapsto\ \  
\frac{1}{|G|} \sum_{g\in G} \rh^*( {}^g(\ta\alpha))
\end{aligned}
$$
induces an isomorphism 
$$(H^{\DOT})^G \stackrel{\sim}{\longrightarrow}\HHD(A).$$
\end{thm}
\vspace{1ex}
\begin{proof}
The statement follows from 
Theorem~\ref{general-S},
Remark~\ref{replace}, 
Proposition~\ref{phi-inverse}, 
and Theorem~\ref{Decomposition}.
\end{proof}

%%%%%%%%%%%%%%%%%%%%%%%%%%%%%
\vspace{3ex}
\begin{remark}
\label{conversion}
We explained explicitly in~\cite{SheplerWitherspoon} how a {\em constant}
Hochschild 2-cocycle defines a graded Hecke algebra.
More generally, Theorem~\ref{S:CGW2} offers a direct conversion
from vector forms in $(H^{\DOT})^G$ to functions on tensor powers of
$S(V)\# G$ representing elements in $\HH^{\DOT}(S(V)\# G)$:
For all $\alpha$ in $H^p$,
$$
  \gamm(\alpha)(f_1\overline{g}_1\ot\cdots\ot f_p\overline{g}_p)=
 \frac{1}{|G|}\sum_{g\in G}\, ^g(\ta\alpha)
(f_1\ot {}^{g_1}f_2\ot\cdots\ot {}^{g_1\cdots g_{p-1}}f_p)\, \overline{g_1\cdots g_p}
$$
for all $f_1,\ldots, f_p$ in $S(V)$ and $g_1,\ldots, g_p$ in $G$.
In particular, if $p=2$, then 
$\gamm(\alpha)(f_1\overline{g}_1\ot f_2\overline{g}_2) = \frac{1}{|G|}
\sum_{g\in G}  {}^g(\ta \alpha) (f_1\ot {}^{g_1}f_2)\,\overline{g_1g_2}$.
\end{remark}
\vspace{3ex}

We now record some inverse isomorphisms that will facilitate finding bracket formulas.
For each $g$ in $G$, let $\{v_1,\ldots, v_n\}$ be a basis
of $V$ consisting of eigenvectors of $g$ and let 
$$\Proj_{H_g}:C_g^{\DOT}\rightarrow H_g^{\DOT}$$ 
be the map which takes any $\CC$-basis element
$ (v_1^{m_1}\cdots v_n^{m_n})\, \overline{g}
\ot v_{i_1}\wedge\cdots\wedge v_{i_p}$
of $C_g^{\DOT}$ to itself if it lies in $H_g^{\DOT}$ and to zero otherwise.
Set $\projH:=\bigoplus_{g\in G} \Proj_{H_g}$.  Then
$$\projH: C^{\DOT}\rightarrow H^{\DOT}\ $$
projects each cocycle in $C^{\DOT}$
to its cohomology class representative in $H^{\DOT}$.
(See, for example, the
computation of~\cite[Section 3.2]{Farinati}.) 
In fact, the map $\projH$ splits with respect
to the canonical embedding $H^{\DOT} \hookrightarrow C^{\DOT}$.
Now consider projection just on the polynomial part of
a form: Let
$$
\Proj_{V^g}:S(V)\overline{g}\ot \Wedge^{\DOT}V^*
\ \rightarrow \
S(V^g)\overline{g}\ot \Wedge^{\DOT}V^*
$$
be the canonical projection arising from the isomorphism
$S(V^g)\cong S(V)/I((V^g)^\perp)$ for each $g$ in $G$.
Then 
$\Proj_{H_g} = \Proj_{H_g} \circ \Proj_{V^g}$ and thus
any cocycle in $C_g$ which vanishes under $\Proj_{V^g}$ is zero in cohomology.
(In fact, we can write $\Proj_H$ as a composition
of $\bigoplus_g \Proj_{V^g}$ with a similar projection map on
the exterior algebra.)

The inclusion map $S(V)\hookrightarrow A$ again induces 
a restriction map (as in Section~\ref{Liftingbracket}):
$$
\res:
\Hom_{\CC}(A^{\otimes\DOT},A)\rightarrow\Hom_{\CC}(S(V)^{\otimes\DOT},A).
$$

Theorem~\ref{general-S2} 
implies that the following compositions 
(for $S=S(V)$) induce inverse isomorphisms
on cohomology,
as asserted in the next theorem:
$$
\begin{xy}
\xymatrix{ 
 \Hom_{\CC}(A^{\ot p},A) \  
\ar@<4pt>[r]^{\res\hspace{.3cm}} &
\Hom_{\CC}(S^{\ot p},A)^G \ 
\ar@{_{(}->}@<5pt>[r] \ 
\ar@<4pt>[l]^{\rh^*} &
\ \Hom_{\CC}(S^{\ot p}, A) \  
\ar@<4pt>[r]^{\hspace{1.1cm}\Phi^*} 
\ar@<4pt>[l]^{\hspace{.3cm}\rey} &
\ C^p \ 
\ar@<4pt>[r]^{\projH} 
\ar@<4pt>[l]^{\hspace{1.1cm}\ta} &
\ H^p 
\ar@{_{(}->}@<5pt>[l] \ .
}
\end{xy}
$$

%%%%%%%%%%%%%%%%%%%%%%%%%%%%%%%%%%%%%%%%%%%%%%%%%%%%%%%%%
\vspace{2ex}
\begin{thm}\label{gamma-prime}
Let $A=S(V)\# G$. 
The maps 
$$\gamm':=\projH\circ\Phi^* \circ \res 
\quad\quad\text{and}\quad\quad
\gamm:=\rh^*\circ \rey \circ \ta$$  
\begin{xy}
\xymatrix{ 
\text{on cochains, }
 \Hom_{\CC}(A^{\ot p},A) \  
\ar@<2pt>[r]^{\hspace{0cm}\Gamma'\hspace{.4cm}} &
\ H^p \ \ar@<2pt>[l]^{\hspace{0cm}\Gamma\hspace{.4cm}}  ,
\text{ induce inverse isomorphisms }
}
\end{xy}
$$\HH^p(A) \stackrel{\sim}{\longleftrightarrow} (H^p)^G.$$
\end{thm}
\vspace{2ex}

The next theorem, a consequence of Theorem~\ref{general-S3}, 
describes the graded Lie bracket on
$(H^{\DOT})^G$ induced by the Gerstenhaber bracket on $\HH^{\DOT}(S(V)\# G)$.

\vspace{2ex}
%%%%%%%%%%%%%%%%%%%%%%%%%%%%%%%%%%%
\begin{thm}\label{bracketunderiso}
The Gerstenhaber bracket on $\HH^{\DOT}(S(V)\# G)$ induces
the following graded Lie bracket 
on $(H^{\DOT})^G$ under the isomorphism
$(H^{\DOT})^G\cong \HH^{\DOT}(S(V)\# G)$ of  Theorem~\ref{gamma-prime}: 
For $\alpha, \beta$ in $(H^{\DOT})^G$,
$$
[\alpha, \beta]
=\frac{1}{|G|^2}\ \projH \sum_{a,b\in G}
\Phi^*  [\rh^*\, ^a(\ta\alpha),  \rh^*\, ^b(\ta\beta)]
$$
where the bracket
$[\ \ , \ \, ]$ on the right is the Gerstenhaber bracket on  $\HH^{\DOT}(S(V)\#G)$
given at the cochain level.
\end{thm}
%%%%%%%%%%%%%%%%%%%%%%%%%%%%%%%%%%%%%%%%%%%%%%%%
\vspace{2ex}

In the remainder of this section, 
we use the above theorem to give formulas for
the Gerstenhaber bracket on $(H^{\DOT})^G\cong \HH^{\DOT}(S(V)\# G)$.  
We introduce notation for a prebracket at the cochain level to aid computations
and allow for an explicit, closed formula.
Recall that we have fixed a basis $B_g$ of $V$ consisting of eigenvectors of
$g$ for each $g$ in $G$.  
These choices are for computational convenience;
our results do not depend on the choices.
For a multi-index $I=(i_1, \ldots, i_m)$, we write
$d v_I$ for $v_{i_1}^*\wedge\cdots\wedge v_{i_m}^*$.
In formulas below, note that we sum over {\em all} multi-indices $I$ of 
a given length, not just those with indices of increasing order.
We use the untagged version $\ttau$ of the map $\ta$ with image in 
$S(V)$ (see Definition~\ref{ttau}) to highlight the group elements
appearing in various formulas.

%%%%%%%%%%%%%%%%%%%%%%%%%%%%%%%%%%%%%%%%%%%%%%%%%%5
\vspace{3ex}
\begin{defn}
\label{prebracketformula}
Define a prebracket (bilinear map) on the set of cochains $C^{\DOT}$
defined in (\ref{C-cochains}),
$$
\begin{aligned}
\hphantom{x}[[\ ,\ ]]%_{\rule[0ex]{0ex}{2ex}B_g, B_h}
: 
           C_g^p \times C_h^q \ \ &\rightarrow \ \ C_{gh}^{p+q-1}+C_{hg}^{p+q-1},
\end{aligned}
$$ 
for $g$ and $h$ in $G$,
depending on a basis of eigenvectors $B_1$ for $g$
and a basis of eigenvectors $B_2$ for $h$ as follows. 
For $\alpha$ in $C_g^p$ and $\beta$ in $C_h^q$, define
$\alpha\overline{\circ}\beta$ to be 
$$
\sum_{\substack{I=(i_1,\ldots,i_m)\vspace{.5ex}\\ 1\leq k\leq p}}
   (-1)^{(q-1)(k-1)}
\ttau_{g, B_1}(\alpha)\bigl(v_{i_1}\otimes\cdots\otimes v_{i_{k-1}}\otimes
f_h^{(i_k)}\otimes\ ^hv_{i_{k+q}} \otimes\cdots\otimes\ ^hv_{i_m}\bigr)
\, \overline{gh}\otimes d v_I,
$$
where $m=p+q-1$, $
f_h^{(i_k)}:= \ttau_{h, B_2}(\beta)
( v_{i_k}\otimes\cdots\otimes v_{i_{k+q-1}}) ,
$
and $v_1,\ldots,v_n$ is any basis of $V$.
Define
$$
\begin{aligned}
\hphantom{x}
[[\alpha, \beta]]_{(\rule[0ex]{0ex}{2ex}B_1, B_2)}
= \alpha\,\overline{\circ}\, \beta 
-(-1)^{(p-1)(q-1)}\,
\beta \,\overline{\circ}\, \alpha\, .
\end{aligned}
$$
\end{defn}
%%%%%%%%%%%%%%%%%%%%%%%%%%%%%%%%%%%%%%%
\vspace{3ex}
%%%%%%%%%%%%%%%%%%%%%%%%%%%
\begin{remark}
\label{caution}
When computing brackets, one may be tempted to seek results 
by working with just 
$\alpha\overline{\circ}\, \beta$ (or just $\beta\,\overline{\circ}\,\alpha$)
and extending by symmetry.  However the 
operation $\overline{\circ}$ is not defined on cohomology.  Furthermore,
one must exercise care in treating the operation $\overline{\circ}$ alone
(e.g., see the proof of Theorem~\ref{bracketalwayszero} below), as one risks 
covertly changing the bases used to apply  $\ttau$
in the middle of a bracket calculation. Similarly, one must
exercise care when examining the bracket of two cocycles summand by summand,
although the bracket is linear:  The bases used to apply $\ttau$
should not depend on the pair of summands considered.
The maps $\ta$ and $\ttau$ are independent
of the choices of bases used when taking brackets of cohomology classes, but a choice
should be made once and for all throughout the whole calculation of a
Gerstenhaber bracket.
\end{remark}
%%%%%%%%%%%%%%%%%%%%%%%%%%%%%%%%%%%%%%%%%%%%%%%%%%%%%%%%%
\vspace{3ex}
We are particularly interested in brackets and prebrackets of elements of cohomological degree~$2$.  
For $\alpha$ in $C_g^2$ and $\beta$ in $C_h^2$, the above definition gives 
\begin{equation}\label{special}
\begin{aligned}
\hphantom{x}
&[[ \alpha, \beta]]_{(\rule[0ex]{0ex}{2ex}B_1, B_2)} = \\
%(1\otimes 1 \otimes v_i & \wedge v_j \wedge v_k)\\
& \hspace{-2ex} \sum_{1\leq i,j,k\leq n}
\Bigl[
 \ttau_{1}{\alpha}\Bigl(\ttau_{2}{\beta}
(v_{i}\otimes v_{j}) \otimes\, ^hv_{k}\Bigr)\, \overline{gh}
\ - \ 
\ttau_{1}{\alpha}\Bigl(v_{i} \otimes
\ttau_{2}{\beta}(v_{j}\otimes v_{k})\Bigr)\, \overline{gh}
\\
& \quad \ \ + 
 \ttau_2{\beta}\Bigl(\ttau_1{\alpha}(v_{i} \otimes 
v_{j})\otimes\, ^gv_{k}\Bigr) \, \overline{hg}
\ - \
\ttau_2{\beta}\Bigl(v_{i}\otimes\ttau_1{\alpha}
(v_{j}\otimes v_{k})\Bigr)\, \overline{hg}\Bigr]
\otimes (v_{i}^*\wedge v_{j}^* \wedge v_{k}^*)  
\end{aligned}
\end{equation}
where $v_1,\ldots, v_n$ is any basis of $V$ and $\ta_1:=\ta_{g, B_1}$
and $\ta_2:=\ta_{h, B_2}$.

%%%%%%%%%%%%%%%%%%%%%%%%%%%%%%%%%%%%%%%%%%%%%%
\begin{remark}\label{bracketindependent}
\vspace{3ex}
The theorem below articulates the Gerstenhaber bracket in terms of 
a fixed basis $B_g$ of $V$ for each $g$ in $G$.
Although each individual prebracket~(\ref{prebracketformula}) depends on
these choices, the formula for the bracket given below
is {\em independent of these choices} (by Remark~\ref{dependsonbasis}).
\end{remark}

%%%%%%%%%%%%%%%%%%%%%%%%%%%%%%%%%%%%%%%%%%%%%%%%%%5
\vspace{2ex}
\begin{theorem}\label{closedGB}
Definition~\ref{prebracketformula} gives a formula
for the Gerstenhaber bracket on $\HH^{\DOT}(S(V)\#G)$ as 
realized on $(H^{\DOT})^G$:
For $g,h$ in $G$ and $\alpha$ in $ H_g^{\DOT}$ and $\beta$ in $H^{\DOT}_h$,
$$
\hphantom{x}
[\alpha, \beta]
=\frac{1}{|G|^2}\  \projH
\sum_{\substack{a,b\in G}}
[[\, ^{a}\alpha,\, ^{b}\beta\, ]]_{(\rule[1ex]{0ex}{1ex}\, ^{a}B_{g}, \,
  ^{b}B_{h})}\ .
$$
\end{theorem}
%%%%%%%%%%%%%%%%%%%%%%%%%%%%%%%%%%%%%%
\vspace{1ex}
\begin{proof}
By Theorem~\ref{bracketunderiso},
the Gerstenhaber bracket $[\ \ , \ \, ]$ 
on $\HH^{\DOT}(S(V)\# G)$ induces the following bracket 
on $H^{\DOT}$:
$$
[\alpha, \beta]
=\frac{1}{|G|^2}\ \projH\sum_{a,b\in G}
\Phi^*[\rh^*\, ^{a}(\ta\alpha),\rh^*\, ^b(\ta\beta)]\ .
$$
We show that this formula is exactly that claimed 
in the special case when $\alpha$ and $\beta$ 
have cohomological degree $2$; the general case follows
analogously.  

Let $v_1, \ldots, v_n$ be any basis of $V$.  
We determine
$
\Phi^* [\rh^*\, ^a(\ta\alpha),  \rh^*\, ^b(\ta\beta)]
$
explicitly in the case $a=b=1_G$
as an element of
$$
\left( S(V)\overline{gh} + S(V)\overline{hg}\, \right)
\ot \Wedge^3 V^*
\cong
\Hom_{\CC}\left(\Wedge^3V,\ S(V)\overline{gh} + S(V)\overline{hg}\, \right)
$$
by evaluating on input of the form $v_i \wedge v_j \wedge v_k$.
The computation for general $a,b$ in $G$ is similar, using (for example)
${}^a(\ta\alpha)
= \ta_{aga^{-1},  {}^aB_g}( {}^a\alpha)$
(see~\cite[Proposition~3.8]{paper1}).
For all $i,j,k$ (see Equation~(\ref{phik})), 
$$
\begin{aligned}
\hphantom{x}
\Phi^* & \left[\rh^*\, (\ta\alpha), \rh^*\, (\ta\beta)\right]\
(v_i \wedge v_j \wedge v_k) \\
&\hspace{10ex} = \ 
\left[\rh^*\, (\ta\alpha),\rh^*\, (\ta\beta)\right]\
\Phi(v_i \wedge v_j \wedge v_k)\\
&\hspace{10ex} = \ 
\sum_{\pi\in \Sym_3}\sgn(\pi)\ 
\left[\rh^*\, (\ta\alpha),\rh^*\, (\ta\beta)\right]
    (v_{\pi(i)}\otimes v_{\pi(j)}\otimes v_{\pi(k)})\ .
\\
\end{aligned}
$$
We expand the Gerstenhaber bracket on $\HH^{\DOT}(S(V)\#G)$ to obtain 
\begin{equation}
\label{}
\begin{aligned}
 \sum_{\pi\in \Sym_3}\sgn(\pi)\ 
\biggl[
&\rh^*(\ta\alpha)
\Bigl(\rh^*(\ta\beta)(v_{\pi(i)}\otimes v_{\pi(j)}) \otimes v_{\pi(k)}\Bigr)\\
\ - \ 
&\rh^*(\ta\alpha)
\Bigl(v_{\pi(i)}\ot \rh^*(\ta\beta)(v_{\pi(j)}\otimes v_{\pi(k)}) \Bigr) \\
\ \ + \
&\rh^*(\ta\beta)
\Bigl(\rh^*(\ta\alpha)(v_{\pi(i)}\otimes v_{\pi(j)}) \otimes v_{\pi(k)}\Bigr)\\
\ - \ 
&\rh^*(\ta\beta)
\Bigl(v_{\pi(i)}\ot \rh^*(\ta\alpha)(v_{\pi(j)}\otimes v_{\pi(k)})
\Bigr) \biggr]\ .
\end{aligned}
\end{equation}

But for any $w_1, w_2, w_3$ in $V$,
$$
\begin{aligned}
\rh^*(\ta \alpha)\Bigl(
\rh^*(\ta\beta)(w_1\ot w_2)\ot w_3\Bigr)
\ =& \
\rh^*\bigl(\ta\alpha\bigr)\Bigl(
\ta\beta(w_1\ot w_2)\ot w_3\Bigr)\\
\ =& \
\ttau\hspace{.2ex}\alpha\biggl(
\ttau\hspace{.2ex}\beta(w_1\ot w_2)\ot\, ^hw_3\biggr)\ \overline{gh}.\\
\end{aligned}
$$
We obtain similar expressions for the other terms arising in the bracket,
and hence the above sum is simply
$$
[[\, \alpha,\, \beta\, ]]_{(\rule[1ex]{0ex}{1ex}B_g, \, B_h)}\ 
(v_i \wedge v_j \wedge v_k)\ .
$$
Similarly, we see that for all $a,b$ in $G$,
$$
\Phi^* \left[\rh^*\, ^a(\ta\alpha),  \rh^*\, ^b(\ta\beta)\right]\
=[[\, ^a\alpha,\, ^b\beta\, ]]_{(\rule[1ex]{0ex}{1ex}\,^aB_g,\,^bB_h)}\ ,
$$
and the result follows.
\end{proof}
%%%%%%%%%%%%%%%%%%%%%%%%%%%%%%%%%%%%%%

\vspace{3ex}

%%%%%%%%%%%%%%%%%%%%%%%%%%%%%%%%%%%%%%%%%%
\begin{remark}\label{extend}
The theorem above actually  gives a formula for a bracket on $H^{\DOT}$, \
not just on
$(H^{\DOT})^G\cong \HHD(S(V)\#G)$.
One can show that this extension to all of $H^{\DOT}$ agrees with
the composition of the Reynolds operator on $H^{\DOT}$
with the Gerstenhaber bracket on $(H^{\DOT})^G$:
By Theorem~\ref{closedGB} and Remark~\ref{RR},
$[\alpha,\beta]$ is $G$-invariant and
cohomologous to $[\rey(\alpha),\rey(\beta)]$
for all $\alpha$, $\beta$ in $H^{\DOT}$.
Hence, the extension to $H^{\DOT}$ is artifical in some sense.
Indeed, a natural Gerstenhaber bracket on all of $H^{\DOT}$ does not make sense,
as this space does not present itself as the Hochschild cohomology
of an algebra with coefficients in that same algebra.  
This idea may be used to explain the formula of Theorem~\ref{closedGB}:
When applying Theorem~\ref{closedGB} to $G$-invariant elements $\alpha$ and $\beta$, 
we might write $\alpha =\sum_{c\in [G/Z(g)]} {}^c\alpha'$
and $\beta = \sum_{d\in [G/Z(h)]} {}^d \beta'$, where $\alpha' \in
(H^{\DOT}_g)^{Z(g)}$ and $\beta' \in (H^{\DOT}_h)^{Z(h)}$ are representative
summands.
(Here, $[G/A]$ is a set of representatives of the cosets $G/A$
for any subgroup $A$ of $G$.)
Then 
$$
\hphantom{x}
[\alpha, \beta]
=\frac{1}{|G|^2}\  \projH
\sum_{\substack{a,b\in G \\ \rule{0ex}{1.5ex}
c\in [G/Z(g)], \ d\in [G/Z(h)]}}
[[\, ^{ac}\alpha',\, ^{bd}\beta'\, ]]_
{(\rule[1ex]{0ex}{1ex}\, ^{a}B_{cgc^{-1}}, \, ^{b}B_{dhd^{-1}})}\ .
$$
However this more complicated expression is cohomologous to that of
Theorem~\ref{closedGB},
since  $\rey\circ\ta(\alpha)$ is cohomologous to $\rey\circ\ta\circ\rey (\alpha)$
(see Remark~\ref{RR}, noting that $\Psi^*=\ta\,$).
\end{remark}
\vspace{3ex}

%%%%%%%%%%%%%%%%%%%%%%%%%%%%%%%%%%%%5
%%%%%%%%%%%%%%%%%%%%%%%%%%%%%%%%%%%%5
\section{Explicit bracket formulas}\label{ebf}

Let $G$ be a finite group and let $V$ be a finite dimensional $\CC G$-module. 
In this section, 
we give a closed formula for the prebracket on $\HH^{\DOT}(S(V)\# G)$.
We shall use this formula in later sections to obtain new results on zero brackets
and to locate noncommutative Poisson structures.
We also recover the classical Schouten-Nijenhuis bracket in this section
and we give an example.

The prebracket of Definition~\ref{prebracketformula} 
simplifies enourmously when we work with a bases of simultaneous
eigenvectors for  $g$ and $h$ in $G$.
Indeed, Remark~\ref{zeroterm} predicts that most terms of Definition~\ref{prebracketformula}
 are zero.  We capitalize on this idea in the next theorem and corollaries.
When the actions of $g$ and $h$ are not simultaneously diagonalizable, we 
enact a change of basis at various points
of the calculation of the prebracket to never-the-less
take advantage of Remark~\ref{zeroterm}. 
The following proof shows how to keep track of the effect
on the map $\ttau$ (mindful of cautionary Remark~\ref{caution}).

First, some notation. For $g,h$ in $G$, let $M=M^{g,h}$ 
be the change of basis matrix between $B_g$ and $B_h$:
For $B_g=\{ w_1,\ldots,w_n\}$
and $B_h=\{v_1,\ldots, v_n\}$ (our fixed bases of eigenvectors of $g$ and
$h$, respectively),
set $M=(a_{ij})$ where, for $i=1,\ldots, n$,
\begin{equation}\label{base-change}
  v_i = a_{1i} w_1+\cdots +a_{ni}w_n.
\end{equation}
The formula below involves  determinants of certain submatrices
of $M$. If $I$ and $J$ are (ordered) subsets of $\{1,\ldots,n\}$, denote by
$M_{I;J}$ the submatrix of $M$ obtained by deleting all rows
except those
indexed by $I$ and deleting all columns except those indexed by $J$.
Recall that
$d v_I:=v_{i_1}^*\wedge\cdots\wedge v_{i_m}^*$
for a multi-index $I=({i_1},\ldots, {i_m})$, 
not necessarily in increasing order.
Below we use the operation $\overline{\circ}$ giving the prebracket
of Definition~\ref{prebracketformula}.
We also use the quantum partial differentiation operators
$\del_{v, \ep}$ (see (\ref{qpd}))
after decomposing $g$ into reflections, $g=s_1\cdots s_n$,
with respect to $B_g$,
i.e., each $s_i$ in $\GL(V)$ is defined by 
$s_i( w_i)={}^g w_i=\ep_i w_i$ and $s_i(w_j)=w_j$ for $j\neq i$.

\vspace{2ex}
%%%%%%%%%%%%%%%%%%%%%%%%%%%%%%%%%%%%%%%%%%%%%%%%%%%%%%%%%%%
%%%%%%%%%%%%%%%%%%%%%%%%%%%%%%%%%%%%%%%%%%%%%%%%%%%%%%%%
\begin{thm}\label{general-circ}
Let $g,h$ lie in $G$ with change of basis matrix $M=M^{g,h}$ as above.
Decompose $g$ into reflections, $g=s_1\cdots s_n$,
with respect to $B_g$ as above.
Let 
$$
\begin{aligned}
\alpha=f_g\overline{g}\ot d w_J\in H^p_g
\quad&\text{where}\quad
J=(j_1<\ldots<j_p)
\quad\text{and}\\
\beta=f_h'\overline{h}\ot d v_L\in H^q_h
\quad&\text{where}\quad
L=(l_1<\ldots< l_q)\ .
\end{aligned}
$$ 
Then $\alpha\overline{\circ}\beta$ is given 
as an element of $C^{p+q-1}_{gh}$ by the following formula:
For $m=p+q-1$ and $I=(i_1<\ldots<i_m)$,
$$
(\alpha\overline{\circ}\, \beta)(v_{i_1}\wedge\cdots\wedge v_{i_m})
= \sum_{1\leq k\leq p}
  (-1)^{\nu(k)}
   \det(M_{J_k;I-L})\
     ^{s_1\cdots s_{j_{k}-1}} (\partial_{j_k}f_h')\ f_g\ \overline{gh},
$$
where $\det(M_{*;I-L})=0$ for $L\not\subset I$, 
$J_k:= (j_1, \ldots,j_{k-1}, j_{k+1}, \ldots, j_p)$,
$\del_{i}=\del_{\ep_i, w_i}$, and 
$\nu(k) = 1-q-k+\lambda(1)+\cdots +\lambda(q) - \binom{q}{2}$
for $\lambda$ defined by $l_s = i_{\lambda(s)}$.
\end{thm}
%%%%%%%%%%%%%%%%%%%%%%%%%%%%%%%%%%%%%%%%%%%%%%%%%%%%%%%%%%%
\vspace{1ex}
\begin{proof}
By Definition~\ref{prebracketformula},
$(\alpha\overline{\circ}\, \beta)(v_{i_1}\wedge\cdots\wedge v_{i_m})$ is equal to
$$\sum_{\stackrel{1\leq k\leq p}{\rule{0ex}{1.5ex}\pi\in \Sym_m}}
  (-1)^{(q-1)(k-1)}(\sgn\pi)\ (\ta_g\alpha)(v_{i_{\pi(1)}}\ot\cdots\ot v_{i_{\pi(k-1)}}\ot
    (f_h')^{(i_{\pi(k)})} \ot  {}^hv_{i_{\pi({k+q})}}\ot\cdots \ot {}^hv_{i_{\pi(m)}})
   \overline{h},
$$
where $(f_h')^{(i_{\pi(k)})}=(\ttau_h\beta)
(v_{i_{\pi(k)}}\ot\cdots\ot v_{i_{\pi(k+q-1)}})$.  
Fix $k$ and note that by Remark~\ref{zeroterm},
$(f_h')^{(i_\pi(k))}$ is nonzero only if $L=(i_{\pi(k)},\ldots,i_{\pi({k+q-1})})$.
Hence, we may restrict the sum to those permutations $\pi$ for which 
$i_{\pi(k)}=l_1,\ldots, i_{\pi({k+q-1})}=l_q$. We identify this set of
permutations
with the symmetric group $\Sym_{p-1}$ (of permutations on the set
$\{1,\ldots,k-1,k+q,\ldots,m\}$) in the standard way.
Under this identification, the factor $\sgn(\pi)$ changes by
$$(-1)^{\lambda(1)-k} (-1)^{\lambda(2)-(k+1)}\cdots 
(-1)^{\lambda(q)-(k+q-1)}\ .$$
For such permutations $\pi$, the vectors $v_{i_{\pi({k+q})}},\ldots, v_{i_{\pi(m)}}$ lie in $V^h$
(as the exterior part of $\beta\in H_h^q$ 
includes a volume form on $(V^h)^\perp$, and hence
$v_i\in V^h$ for $i\notin L$).
Therefore $ {}^hv_{i_{\pi(k+q)}} = v_{i_{\pi(k+q)}},\ldots, {}^hv_{i_{\pi(m)}}=v_{i_{\pi(m)}}$.
Note also that $(f_h')^{(i_{\pi(k)})}=f_h'$ for all such $\pi$.
We now invoke the change of basis
$$
   v_{i_{\pi(s)}}=a_{1,i_{\pi(s)}}w_1+\cdots +a_{n,i_{\pi(s)}}w_n
$$
for $s=1,\ldots,k-1, k+q,\ldots, m$ to obtain\\
\hphantom{xxx}
$ \ta_g(\alpha)(v_{i_{\pi(1)}}\ot\cdots\ot v_{i_{\pi({k-1})}}\ot
    f_h' \ot v_{i_{\pi({k+q})}}\ot\cdots\ot  v_{i_{\pi(m)}})$
$$\hphantom{xx}= a_{j_1,i_{\pi(1)}}\cdots
   a_{j_{k-1},i_{\pi({k-1})}}\cdot a_{j_{k+1},i_{\pi({k+q})}}\cdots 
   a_{j_p,i_{\pi(m)}}\  ^{s_1\cdots s_{j_{k}-1}}(\partial_{j_k}
   f_h') f_g\,\overline{g}.
$$

Thus the $k$-th summand of 
  $(\alpha\overline{\circ}\,\beta)(v_{i_1}\wedge\cdots\wedge v_{i_m})$ 
is a sum over
$\pi\in \Sym_{p-1}$
of $$(-1)^{(q-1)(k-1)+ \lambda(1)+\cdots + \lambda(q)
   -k - (k+1) - \cdots - (k+q-1)} (\sgn\pi)$$ times
$$
   a_{j_1,i_{\pi(1)}}
  \cdots a_{j_{k-1},i_{\pi(i_{k-1})}} a_{j_{k+1},i_{\pi({k+q})}}\cdots a_{j_p,
    i_{\pi(m)}}\ ^{s_1\cdots s_{j_{k}-1}} (\partial_{j_k}f_h')
   f_g \, \overline{gh}.
$$
By definition of  the  determinant  and some sign simplifications,
this is precisely the formula claimed.
\end{proof}
\vspace{1ex}

In the remainder of this section, we refine the formula of Theorem~\ref{general-circ}
in the special case that the actions of $g$ and $h$ commute.
In this case, we may choose
$B_g=B_h$ to be a simultaneous basis of eigenvectors for both $g$ and $h$, 
and the change of basis matrix $M$ in Theorem~\ref{general-circ}
is simply the identity matrix.  
We introduce some notation to capture the single summand that
remains, for each $k$, after the prebracket formula in Definiton~\ref{prebracketformula} 
collapses.
For any two multi-indices $J=(j_1,\ldots,j_p)$ and
$L=(l_1,\ldots,l_q)$, and $k\leq p$, 
set 
$$
\begin{aligned}
I_k :=(j_1,\ldots, j_{k-1},l_1,\ldots, l_q, j_{k+1},\ldots,j_p)\\
I'_k :=(l_1,\ldots, l_{k-1},j_1,\ldots, j_p, l_{k+1},\ldots,l_q).
\end{aligned}
$$

\vspace{2ex}
%%%%%%%%%%%%%%%%%%%%%%%%%%%%%%%%%%%%%%
\begin{cor}\label{prebracket1}
Suppose the actions of $g,h$ in $G$ on $V$ commute and 
$B_g=B_h=\{v_1,\ldots,v_n\}$ is a simultaneous basis of eigenvectors for $g$ and $h$.
Let $\alpha=f_g \overline{g}\otimes d v_J$ lie in $H^p_g$
and $\beta=f_h' \overline{h} \otimes d v_L$ lie in $H^q_h$,
for multi-indices $J=(j_1<\ldots< j_p)$ and $L=(l_1<\ldots< l_q)$.
Then
$$
\begin{aligned}
\hphantom{x}
[[\alpha, \beta]]_{(\rule[1ex]{0ex}{1ex}B_g, B_h)}\quad\quad\quad 
=&\sum_{1\leq k\leq p}
 (-1)^{(q-1)(k-1)} {}^{s_{1}\cdots s_{j_{k}-1}}(\del _{{j_k}}f_h')\ 
f_g\ \overline{gh}\otimes d v_{I_k}
\\
-(-1)^{(p-1)(q-1)}\cdot
 &\sum_{1\leq k\leq q}
  \, (-1)^{(p-1)(k-1)}\ ^{s'_{1}\cdots s'_{l_{k}-1}}(\del '_{{l_k}}f_g)\
f_h'\ \overline{hg}\otimes d v_{I'_k}\\
\end{aligned}
$$
where $s_i$, $s_i'$ are
diagonal reflections (or identity maps) with $^{s_i}v_i=\, ^{g}v_i=\ep_i v_i$
and $^{s'_i}v_i=\, ^hv_i=\ep'_i v_i$, and where
$\del_{i}=\del_{v_i,\ep_i}$  and 
$\del'_{i}=\del_{v_i,\ep_i'}$ .
\end{cor}
\vspace{2ex}
%%%%%%%%%%%%%%%%%%%%%%%%%%%%%%%%%%%%%%%%%%%

If the group elements $g$ and $h$ act as the identity, we
recover the classical Schouten-Nijenhuis bracket, as detailed in the next two corollaries.
Indeed, we view Corollary~\ref{prebracket1} 
as giving a
quantum version of the Schouten-Nijenhuis bracket.
\vspace{2ex}
%%%%%%%%%%%%%%%%%%%%%%%%%%%%%%%%%%%%%%
\begin{cor}\label{prebracket3}
Let $g$ and $h$ in $G$ both act trivially on $V$
(e.g., $g=h=1_G$). 
Let $\alpha=f_g \overline{g}\otimes d v_J$ lie in $H^p_{g}$
and $\beta=f_h' \overline{h} \otimes d v_L$ lie in $H^q_{h}$,
for multi-indices $J=(j_1<\ldots< j_p)$ and $L=(l_1<\ldots< l_q)$.
If $B_g=\{v_1,\ldots,v_n\}=B_h$, then
$$
\begin{aligned}
\hphantom{x}
[[\alpha, \beta]]_{(\rule[1ex]{0ex}{1ex}B_g, B_h)}
=
  & \sum_{1\leq k\leq p}
\ (-1)^{(q-1)(k-1)}
\frac{\del}{\del v_{{j_k}}}(f_h')\ f_g\ \overline{gh}\otimes d v_{I_k}\\
& -(-1)^{(p-1)(q-1)}
 \sum_{1\leq k\leq q}
  \ (-1)^{(p-1)(k-1)}
\frac{\del}{\del v_{{l_k}}}(f_g)\ f_h'\ \overline{hg} \otimes d v_{I'_k}\ .
\end{aligned}
$$
\end{cor}
%%%%%%%%%%%%%%%%%%%%%%%%%%%%%%%%%%%%%%%%%%%%%%%%%%%%%%%%
\vspace{2ex}

If the group $G$ is trivial, then the bracket agrees with the prebracket:

\vspace{2ex}
%%%%%%%%%%%%%%%%%%%%%%%%%5
\begin{cor}\label{sch-nij}
If $G=\{1\}$, we recover the classical Schouten-Nijenhuis bracket:
Let $\alpha = f_1\ot dv_J$ in degree $p$ and $\beta=f_2\ot dv_L$ in degree $q$.
Then
$$
\begin{aligned}
\hphantom{x}
  [\alpha,\beta] =  & \sum_{1\leq k\leq p} (-1)^{(q-1)(k-1)}
   \frac{\partial}{\partial v_{j_k}} (f_2) f_1 \ot dv_{I_k}\\
     & - 
   (-1)^{(p-1)(q-1)} \sum_{1\leq k\leq q} (-1)^{(p-1)(k-1)}\frac{\partial}{\partial v_{l_k}}
  (f_1) f_2 \ot dv_{I_k'}\ .
\end{aligned}
$$ 
\end{cor}
%%%%%%%%%%%%%%%%%%%%%%%%%%%%%%%%%%%%%%%%%%%%%%%%%%%%%%
\vspace{2ex}

%%%%%%%%%%%%%%%%%%%%%%%%%%%%%%%%%%%%%%%%%%%%%%%%%%%%%%%%%%%%%%%%%%%%%%%%%%%%%%%%%
%%%%%%%%%%%%%%%%%%%%%%%%%%%%%%%%%%%%%%%%%%%%%%%%%%%%%%%%%%%%%%%%%%%%%%%%%%%%%%%%%
%%%%%%%%%%%%%%%%%%%%%%%%%%%%%%%%%%%%%%%%%%%%%%%%%%%%%%%%%%%%%%%%%%%%%%%%%%%%%%%%%
%%%%%%%%%%%%%%%%%%%%%%%%%%%%%%%%%%%%%%%%%%%%%%%%%%%%%%%%%%%%%%%%%%%%%%%%%%%%%%%%%

For abelian groups, brackets enjoy a combinatorial description:

\vspace{3ex}
\begin{example}\label{abelian}
Let $G$ be an abelian group acting on $V$ with $\dim V\geq 3$.
Let $B=\{v_1,\ldots,v_n\}$, a simultaneous basis of eigenvectors for $G$.
For $i\in\{1,\ldots,n\}$, let $\chi_i$ be the character
defined by ${}^{g}v_i = \chi_i(g)\, v_i$  for each $g$ in $G$.
Fix $g,h\in G$ and 
$$\begin{aligned}  
\alpha & =
(v_1^{c_1}v_2^{c_2}v_3^{c_3})\,\overline{g}\ot v_1^*\wedge v_2^*,\\
\beta &= (v_1^{d_1}v_2^{d_2}v_3^{d_3})\,\overline{h}\ot v_2^*\wedge v_3^*\ ,
\end{aligned}
$$
elements of $H^2_g$, $H^2_h$, respectively.
(Note that if $g$ acts nontrivially on $V$, then $c_1=c_2=0$, and if
$h$ acts nontrivially on $V$, then $d_2=d_3=0$.)
A calculation using Theorems~\ref{closedGB} and \ref{general-circ}
gives the bracket:
$$
[\alpha,\beta] \ = \kappa\ \projH\
v_1^{c_1+d_1}v_2^{c_2+d_2 -1} v_3^{c_3+d_3} 
    \overline{gh}\ot v_1^*\wedge
     v_2^*\wedge v_3^* \ ,
$$
where
$$
\kappa\ =\   \hphantom{\, - \ }
\langle \chi_1^{c_1-1}\chi_2^{c_2-1} ,  \chi_3^{-c_3}\rangle
    \langle \chi_1^{d_1},  \chi_2^{1-d_2}\chi_3^{1-d_3}\rangle 
   \ \bigl([c_2]_\ep\ \chi_1(h)^{c_1}
    -  
   \ [d_2]_{\ep'}\ \, \chi_1(g)^{d_1}\bigr) \ .
$$
Here, $\ep={\chi_2(h)}, \ep'={\chi_2(g)}$, 
$[m]_{\lambda}$ is the quantum integer
$1+\lambda+\ldots+\lambda^{m-1}$ (or zero when $m=0$),
and $\langle \hspace{.1in}, \hspace{.1cm} \rangle$ denotes the inner product
of characters on $G$.
The map $\projH$ in particular
projects the polynomial coefficient onto $S(V^{gh})$:
$v_i\mapsto  v_i$ if $\chi_i(gh)=1$,
$v_i \mapsto 0$ otherwise.
Thus $[\alpha,\beta]$ is usually $0$, but can be nonzero 
as a consequence of the orthogonality relations on characters of finite groups
(see Proposition~\ref{nonzeroabelianbracket} below).
\end{example}
\vspace{3ex}

%%%%%%%%%%%%%%%%%%%%%%%%%%%%%%%%%%%%%%%%%%%%%%
%%%%%%%%%%%%%%%%%%%%%%%%%%%%%%%%%%%%%%%%%%%%%%%%%%%5
%%%%%%%%%%%%%%%%%%%%%%%%%%%%%%%%%%%%%%%%%%%%%%%%%%%%%%%
\section{Zero brackets}\label{sec:zero}

Let $G$ be a finite group and let $V$ be a finite dimensional $\CC G$-module. 
Every deformation arises from a Hochschild $2$-cocycle $\alpha$ whose
square bracket $[\alpha,\alpha]$ is a coboundary.  
We use the formulas of Section~\ref{ebf} to 
now determine some conditions under which 
brackets are zero.  In the process, we take advantage of
our depiction of cohomology automorphisms 
in terms of {\em quantum partial derivatives}.
We begin with an easy corollary of our formulation of 
the Gerstenhaber bracket as a sum of prebrackets.

We say that a cochain $\alpha$ in $C^{\DOT}$ is {\bf constant}
if the polynomial coefficient of $\alpha$ is constant, i.e.,
if $\alpha$ lies in the subspace
$\bigoplus_{g\in G} \CC\, \overline{g} \ot \Wedge V^*$.
We showed in~\cite{SheplerWitherspoon} 
that any constant Hochschild $2$-cocycle lifts to a graded Hecke algebra.
As graded Hecke algebras are deformations of
$S(V)\#G$ (see Section~\ref{previouspaper}),
every constant Hochschild $2$-cocycle
defines a noncommutative Poisson structure.
The following consequence of Theorem~\ref{closedGB}
extends this result to {\em arbitrary} cohomological degree.

%%%%%%%%%%%%%%%%%%%%
\vspace{2ex}
\begin{thm}\label{polydeg0}
Suppose $\alpha$ and $\beta$ in $H^{\DOT}$ are constant.
Then $[\alpha,\beta]=0$.
\end{thm}
\vspace{1ex}
\begin{proof}
Theorem~\ref{closedGB} gives the Gerstenhaber bracket
as a sum of prebrackets of cochains 
$^a\alpha$ and $^b\beta$ for $a,b$ in $G$.  
Definition~\ref{prebracketformula} expresses
the prebracket in terms of quantum partial differentiation of 
polynomial coefficients of cochains (see Definitions~\ref{tau}
and \ref{ttau}).  
As the cochains $^a\alpha$ and $^b\beta$ are constant, 
any partial derivative of their coefficients 
is zero.  Hence, each prebracket is zero.
\end{proof}
\vspace{1ex}

As an immediate consequence of the theorem, we obtain:

\vspace{2ex}
\begin{cor}
Any {\em constant} Hochschild cocycle in $\HHD(S(V)\#G)$
defines a noncommutative Poisson structure. 
\end{cor}
\vspace{2ex}

We need a quick linear algebra lemma in order to give more results
on zero brackets.

%%%%%%%%%%%%%%%%%%%%%%%%%%%%%%%%%%%%%%%%%%%%%%%
\vspace{2ex}
\begin{lemma}\label{zerodet}
Let $v_1,\ldots,v_n$ and $w_1,\ldots, w_n$ be two
bases of $V$ with change of basis matrix $M$ determined by Equation~\ref{base-change}.
If $J$ and $L$ are two multi-indices with
$$
 \Span_{\CC}\{ w_{j}\, |\, j \in J \}
   \cap \Span_{\CC} \{ v_{l}\,|\, l \in L\}\ \neq \{0\}\ ,
$$
then $\det (M_{J; L'}) =0$
for any $L'$ with $L\cap L'=\emptyset$.
\end{lemma}
\vspace{1ex}
\begin{proof}
Write 
$
   \sum_{j\in J} c_jw_j = \sum_{l\in L} d_lv_l
$
for some scalars $c_j$, $d_l$, not all 0.
We substitute $v_i=a_{1i}w_i+\ldots+a_{ni}w_n$ 
to see that $\sum_{l\in L}d_la_{il}=0$ for $i\notin J$.
If we delete the rows of $M$ indexed by $J$, these
equations give a linear dependence relation among
the columns indexed by $L$.
Thus $\det(M_{J;L'})=0$. 
\end{proof}
\vspace{1ex}

As a consequence, we have the following proposition.

\vspace{2ex}
%%%%%%%%%%%%%%%%%%%%%%%%%%%%%%%%%%%%%%%%%%%%%%%%%%%%%%%%%%%%%%%%%%%%%%%%%%%%%%%%%%%
\begin{prop}
Suppose $g,h$ lie in $G$ and $(V^g)^\perp\cap (V^h)^\perp$ is nontrivial and fixed
setwise by $G$.  Let $\alpha$ and $\beta$ be elements of
$(H^{\DOT})^G$ supported on the conjugacy classes of $g$ and $h$,
respectively.  Then
$$[\alpha,\beta]=0\ .$$
\end{prop}
\vspace{1ex}
%%%%%%%%%%%%%%%%%%%%%%%%%%%%%%%%%%%%%
\begin{proof}
Let $\alpha'$ be the summand of $\alpha$ supported on $g$ itself and
let $\beta'$ be the summand of $\beta$ supported on $h$.
Then $|Z(g)|\, \alpha=\rey\alpha'$.
(Indeed, since $\alpha$ is $G$-invariant,
$\alpha'$ is $Z(g)$-invariant and
$\alpha=\sum_{c} \,^c\alpha'$, a sum over coset representatives $c$
of $G/Z(g)$.) 
Similarly, $|Z(h)|\,\beta=\rey\beta'$.
Remark~\ref{RR} then implies that
$$
(|Z(g)|\,|Z(h)|)\,[\alpha,\beta] = [\rey\alpha',\rey\beta'\,]=[\alpha',\beta'\,]\ ,
$$
and hence we compute $[\alpha',\beta'\,]$.
Set $W:=(V^g)^\perp\cap (V^h)^\perp$
and recall that the Hermitian form on $V$ is $G$-invariant.
By Remark~\ref{bracketindependent}, we may assume
that
$B_g=\{w_1, w_2, \ldots, w_n\}$ and
$B_h=\{v_1, v_2, \ldots, v_n\}$
are orthogonal bases with $w_1, v_1$ in $W$.
Without loss of generality, suppose
$$
\alpha'=f_g \overline{g}\ot w_{j_1}^*\wedge\cdots\wedge w_{j_p}^*
\quad\text{and}\quad
 \beta'=f_h' \overline{h}\ot v_{l_1}^*\wedge\cdots\wedge v_{l_q}^*
$$
in $H^{\DOT}$ where $f_g\in S(V^g)$, $f_h'\in S(V^h)$
and $j_1=l_1=1$.
As $\beta'$ lies in $H^q_h$, the linear span of $v_{l_1}, \ldots, v_{l_q}$
contains $(V^h)^\perp$ and thus $w_1$. Hence, for $k\neq 1$,
$$
   \Span_{\CC}\{ \,w_{j_i}\, | \, i=1,\ldots, p\,;\ i \neq k\}
   \cap \Span_{\CC} \{ \,v_{l_i}\ | \ i=1,\ldots, q\}
\neq \{0\}\ .
$$
Theorem~\ref{general-circ} and Lemma~\ref{zerodet}
then imply that the $k$-th summand of 
$\alpha'\,\overline{\circ}\,\beta'$ in the prebracket
$[[\,  \alpha', \,  \beta'\,]]_{(\rule[0ex]{0ex}{2ex}\, B_g, \, B_h)}$
is zero for $k\neq 1$ (as $\det M_{J_k; I-L}=0$).  
But the first ($k=1$) summand is also zero,
since $\, w_1=w_{j_1}$ lies in $(V^h)^{\perp}$
while $f_h'$ lies in $V^h$, forcing
the partial derivative of $f_h'$ with respect to $w_1$ and $B_g$ to be zero:
$\del_{1}(f_h')=0$.
(Indeed, since $w_1\in (V^h)^\perp$ and $B_g$ is orthogonal,
$V^h\subset\Span_{\CC}\{w_2,\ldots,w_n\}$:
Were some $w_1+a_2w_2+\ldots a_nw_n$ to lie in $V^h$ (with $a_i\in\CC$),
the inner product
$\left<w_1,w_1\right>=\left<w_1,w_1+a_2w_2+\ldots+a_nw_n\right>$ 
would be zero.)
Hence
$\alpha'\,\overline{\circ}\,\beta'=0$.
By a symmetric argument, interchanging
the role of $g$ and $h$, 
$\beta'\,\overline{\circ}\,\alpha'$ is also zero.

The same arguments apply when we compute the prebracket
$[[\,  ^a\alpha', \,  ^b\beta'\,]]$
with respect to the pair of bases
$(\, ^aB_g, \, ^bB_h)$
for $a,b$ in $G$, as $^a w_1$ and $^b v_1$ both lie in
$W= {}^aW = {}^bW$ by hypothesis and 
\begin{equation}\label{intersection}
   W \subset\ ^a((V^g)^{\perp}) \cap\ ^b((V^h)^{\perp})
   = (V^{aga^{-1}})^{\perp} \cap\ (V^{bhb^{-1}})^{\perp} \ .
\end{equation}
Hence, Theorem~\ref{closedGB} implies that $0=[\alpha', \beta']=[\alpha,\beta]$.
\end{proof}
%%%%%%%%%%%%%%%%%%%%%%%%%%%%%%%%%%%%%%%%%%%%%%%%%%%%%%%%%
\vspace{1ex}

We immediately obtain some interesting corollaries 
for square brackets:
%%%%%%%%%%%%%%%%%%%%%%%%%%%%%%%%%%%%%%%%%%%%%%%%%%%
\vspace{2ex}
\begin{cor}
Let $g$ lie in $G$.
If $V^g\neq V$ is fixed setwise by $G$, then 
any Hochschild cocycle $\alpha$ in $\HH^{\DOT}(S(V)\#G)$ supported on  
conjugates of $g$ defines a noncommutative Poisson bracket:
$$[\alpha, \alpha]=0\ .$$
\end{cor}
\vspace{1ex}
%%%%%%%%%%%%%%%%%%%%%%%%%%%%%%%%%%%%%%%%%%%%%
Note that the hypothesis of the corollary is automatically satisfied for any element
$g$  in the {\em center} of $G$ acting nontrivially on $V$. Thus the square
bracket of any Hochschild cocycle in $\HH^{\DOT}(S(V)\#G)$ supported on
conjugates of such an element $g$ is zero.
%%%%%%%%%%%%%%%%%%%%%%%%%%%%%%%%%%%%%%%%%%%%%%%%%%%

More generally, we have a corollary for brackets of possibly different cocycles:

%%%%%%%%%%%%%%%%%%%%%%%%%%%%%%%%%%%%%%%%%%%%%%%%%%%
\vspace{2ex}
\begin{cor}
Suppose $g,h$ lie in the center of $G$ with
$(V^g)^\perp\cap (V^h)^\perp$ nontrivial.  Then the bracket of
any two Hochschild cocycles $\alpha, \beta$ in $\HH^{\DOT}(S(V)\#G)$ supported on 
conjugates of $g,h$, respectively, is zero:
$$[\alpha,\beta]=0\ .$$
In particular, if $G$ is abelian, such a bracket is always 0 whenever
$(V^g)^\perp\cap (V^h)^\perp$ is nontrivial.
\end{cor}

%%%%%%%%%%%%%%%%%%%%%%%%%%%%%%%%%%%%%%%%%%%%%%%%%%%%%%%%%%%%%5
%%%%%%%%%%%%%%%%%%%%%%%%%%%%%%%%%%%%%%%%%%%%%%%%%%%%%%%%%%%%%5
%%%%%%%%%%%%%%%%%%%%%%%%%%%%%%%%%%%%%%%%%%%%%%%%%%%%%%%%%%%%%5
%%%%%%%%%%%%%%%%%%%%%%%%%%%%%%%%%%%%%%%%%%%%%%%%%%%%%%%%%%%%55
\section{Zero brackets in cohomological degree 2}\label{sec:degree-2}

Let $G$ be a finite group and let $V$ be a finite dimensional $\CC G$-module. 
We are particularly interested in the Gerstenhaber bracket
in cohomological degree $2$, since every deformation of $S(V)\# G$
arises from a Hochschild 2-cocycle $\mu$
with square bracket $[\mu, \mu]$
zero in cohomology, i.e., from a 
noncommutative Poisson structure (see Section~\ref{Prelim}). 
Graded Hecke algebras are deformations of $S(V)\# G$ that
arise from noncommutative Poisson structures of a particular form (see Section~\ref{previouspaper}).
In Theorem~\ref{bracketalwayszero} below, we prove that Hochschild 2-cocycles
supported on group elements acting nontrivially
on $V$ always have zero bracket.  As an immediate consequence, 
if $G$ acts faithfully, then any single Hochschild 2-cocycle supported on nonidentity
group elements defines a noncommutative Poisson structure.  

The proof of Theorem~\ref{bracketalwayszero} rests on the following combinatorial
lemma about Demazure operators and reflections.  
In fact, by enunciating our automorphisms of cohomology
in terms of Demazure operators, we show in this lemma 
that the operation $\overline{\circ}$ 
(see Definitions~\ref{prebracketformula} and~\ref{circ}) 
is essentially zero on 
$2$-cocycles after projecting to cohomology classes under
the hypotheses of the theorem.

\vspace{2ex}
%%%%%%%%%%%%%%%%%%%%%%%%%%%%%%%%%%%%%%%%%%
\begin{lemma}\label{udivides}
Suppose $g$ in $G$ acts as a bireflection, i.e., $\codim V^g =2$. 
For all $w$ in $V$, $\alpha$ in $H^2_g$,
and  $m\geq 0$, the difference $w- {}^gw$ divides 
$$\ttau\hspace{.2ex}\alpha(w^m \ot w )
-\ttau\hspace{.2ex}\alpha(w\ot w^m).$$
\end{lemma}
%%%%%%%%%%%%%%%%%%%%%%%%%%%%%%%%%%%%%5
\vspace{1ex}
\begin{proof}
Let $B_g=\{v_1,\ldots, v_n\}$.
Without loss of generality,
assume $v_1,v_2$ span $(V^g)^\perp$.
Let $\epsilon_1, \epsilon_2$ be the corresponding (nontrivial) eigenvalues of $g$.
Write $\alpha=f_g\overline{g}\otimes v_1^*\wedge v_2^*$ for 
some $f_g$ in $S(V)$.
Let $s_1$ and $s_2$ in $\GL(V)$ be diagonal reflections decomposing $g$, i.e.,
$^{s_1}v_1=\, ^gv_1=\ep_1v_1$, $^{s_1}v_2=v_2$, $^{s_2}v_1=v_1$, 
$^{s_2}v_2=\, ^g(v_2)=\ep_2v_2$,
and $g=s_1s_2=s_2 s_1$.  By Definitions~\ref{tau} and~\ref{ttau},
$\ttau\hspace{.2ex}\alpha(w^m\ot w)-\ttau\hspace{.2ex}\alpha(w\ot w^m)$ 
is the following multiple of $f_g $:
$$
\begin{aligned}
(\partial_1w^m)\ ^{s_1}(\partial_2w) - & (\partial_1w)\ ^{s_1}(\partial_2w^m)\\
= \ 
&
\left(\displaystyle{\frac{w^m-\, ^{s_1}w^m}{v_1-\, ^{s_1}v_1}}\right) 
\ \raisebox{0ex}[3ex]{}^{s_1}\hspace{-1ex}
\left(\frac{w-\, ^{s_2}w}{v_2-\, ^{s_2}v_2}\right) - 
   \left(\frac{w-\, ^{s_1}w}{v_1-\, ^{s_1}v_1}\right) 
\ \raisebox{0ex}[3ex]{}^{s_1}\hspace{-1ex}
\left(\frac{w^m-\, ^{s_2}w^m}{v_2-\,^{s_2}v_2}\right).
\end{aligned}
$$
Since $v_1-\, ^{s_1}v_1=(1-\epsilon_1)v_1$ and $v_2-\, ^{s_2}v_2
=(1-\epsilon_2)v_2$,
we may factor out the scalar $(1-\epsilon_1)^{-1}(1-\epsilon_2)^{-1}$
from each summand, leaving us with 
${\displaystyle\frac{1}{v_1 v_2}}$ times
$$
\begin{aligned}
&\left({w^m-\, ^{s_1}w^m}\right)
\left({\, ^{s_1}w-\, ^{g}w}\right) 
-
\left({w-\, ^{s_1}w}\right)
\left({\, ^{s_1}w^m-\,^{g}w^m}\right)\\
&\ =\
w^m(\,^{s_1}w) - w^m(\,^gw) + (\,^{s_1}w^m)(\,^gw)
   -w(\,^{s_1}w^m) + w(\,^gw^m) -(\,^{s_1}w)(\,^{g}w^m) \\
&\ =\
(\,^{s_1}w)(w^m-\,^gw^m) 
+ (\,^{s_1}w^m)(\,^gw-w)
-w^m(\,^gw) + w^{m+1} -w^{m+1} + w(\,^gw^m)\\
&\ = \ 
(w-\,^gw)(w^m-\,^{s_1}w^m)-(w-\,^{s_1}w)(w^m-\,^gw^m)\\
&\ = \  
(w-\,^gw)(w-\,^{s_1}w)[w^{m-1}+w^{m-2}(\,^{s_1}w)+\ldots+\,^{s_1}w^{m-1}]\\
&\hspace{11ex}
-(w-\,^{s_1}w)(w-\,^gw)[w^{m-1}+w^{m-2}(\,^gw)+\cdots+\,^gw^{m-1}]\\
&\ =\
(w-\,^gw)(w-\,^{s_1}w)\\
&\hspace{10ex}
\cdot \bigl(
w^{m-2}(\,^{s_1}w-\,^gw)+ 
w^{m-3}(\,^{s_1}w^2-\,^gw^2)+\ldots+(\,^{s_1}w^{m-1}-\,^gw^{m-1})
\bigr)\, .\\
\end{aligned}
$$
Now $w-\, ^{s_1}w$ lies in ${\rm im} (1-s_1)=(V^{s_1})^{\perp}={\rm Span}_{\mathbb C}
\{ v_1\}$.  Similarly, for each $i$, $\,^{s_1}w^i-\,^g w^i=
(\,^{s_1}w^i)-\,^{s_2}(\,^{s_1}w^i)$ lies in ${\rm im}(1-s_2)=\Span_{\CC}\{v_2\}$.
Hence, $v_1 v_2$ divides the above expression in $S(V)$ and the quotient
by $v_1v_2$ is divisible by $w-\, ^gw$.  
\end{proof}
\vspace{1ex}
%%%%%%%%%%%%%%%%%%%%%%%%%%%%%%%%%%%%%%%%%%%%%%%%%%5

We are now ready to show that the Gerstenhaber bracket 
is zero on cocycles supported on group elements acting nontrivially.
In the next section, we show that the converse of this theorem
is false and that the hypothesis in the theorem
can not be easily weakened.
Let $K$ denote the kernel of the action of $G$ on $V$.

%%%%%%%%%%%%%%%%%%%%%%%%%%%%%%%%%%%%%%%%%%%%
%%%%%%%%%%%%%%%%%%%%%%%%%%%%%%%%%%%%%%
\vspace{2ex}
\begin{thm}
\label{bracketalwayszero}
The  bracket of any two elements $\alpha$, $\beta$ 
in $\HH^2(S(V)\#G)$ supported off $K$ is zero:
$$
[\alpha,\beta]=0\ .
$$ 
\end{thm}
%%%%%%%%%%%%%%%%%%%%%%%%%%%%%%%%%%%%%
\vspace{1ex}
\begin{proof}
For each $k$ in $G$,
we fix a basis $B_k$ of $V$ consisting of
eigenvectors of $k$ 
so that the first $\codim(V^k)$ vectors in that list span $(V^k)^\perp$.
Then for all $a$ in $G$, the first $\codim(V^{aka^{-1}})$
vectors in the basis $\,^aB_k$  span $(V^{aka^{-1}})^\perp$ as well.

Since the bracket is linear, we may assume without loss of generality
that $\alpha$ and $\beta$ are each supported on the conjugacy class of a single element 
in $G$. In fact, it is enough to consider single summands:
Assume that $\alpha$ lies in $H^2_g$ and $\beta$ lies in $H^2_h$
for some $g$ and $h$ in $G$.
By Theorem~\ref{closedGB},
our Definition~\ref{prebracketformula} gives a closed formula
for the Gerstenhaber bracket on $\HH^{\DOT}(S(V)\#G)$ as realized on $H^{\DOT}$:
$$
\hphantom{x}
[\alpha, \beta]
=\frac{1}{|G|^2}\  \projH
\sum_{\substack{a,b\in G }}
[[\, ^{a}\alpha,\, ^{b}\beta\, ]]_{\rule[1ex]{0ex}{1ex}(\,^{a}B_g,\,^{b} B_h)}\ .
$$
(By Remark~\ref{bracketindependent},
this bracket does not depend on our choices of bases $B_k$ for each $k$ in $G$.)
Write $B_1:=\{w_1,\ldots, w_n\}$ for ${}^aB_g$ 
and $B_2:=\{v_1, \ldots, v_n\}$ for ${}^bB_h$.
We assume first that $a=b=1_G$.
Consider the prebracket $[[\alpha,\, \beta]]:=[[\alpha,\,
    \beta]]_{(\rule[0ex]{0ex}{2ex}B_1, B_2)}$,
a cochain in $C^3_{gh} + C^3_{hg}$,
and suppose that $H^3_{gh}$ or $H^3_{hg}$ is nonzero.
We show that this prebracket is either identically zero 
or projects to zero under the map $\projH$.  

For any $i<j<k$ and
$\ttau_1=\ttau_{B_1,g}$ and $\ttau_2=\ttau_{B_2,h}$,
Equation~(\ref{special}) gives
\begin{equation}
\label{four2}
\begin{aligned}
\hphantom{x}
[[ &\alpha, \beta]]( v_i \wedge v_j \wedge v_k) 
 =  \sum_{\pi\in \Sym_3}\sgn(\pi)\cdot\,\\
&\hspace{3ex} 
\Biggl( 
\ttau_1{\alpha}\Bigl(\ttau_2{\beta}
(v_{\pi(i)}\otimes v_{\pi(j)}) \otimes\, ^h v_{\pi(k)}\Bigr) 
\, \overline{gh}
\, - \, 
\ttau_1{\alpha}\Bigl(v_{\pi(i)} \otimes
\ttau_2{\beta}(v_{\pi(j)}\otimes v_{\pi(k)})\Bigr)
\, \overline{gh}
\\
&\hspace{2ex}+  
\ttau_2{\beta}\Bigl(\ttau_1{\alpha}(v_{\pi(i)} \otimes 
v_{\pi(j)})\otimes\, ^g v_{\pi(k)}\Bigr) 
\,\overline{hg} - \,
\ttau_2{\beta}\Bigl(v_{\pi(i)}\otimes\ttau_1{\alpha}
(v_{\pi(j)}\otimes v_{\pi(k)})\Bigr)
\, \overline{hg}\ \Biggr)\ .
\end{aligned}
\end{equation}
We may assume $g$ and $h$ act as bireflections, i.e.,
$\codim V^g=\codim V^h=2$,
else $H^2_g$ or $H^2_h$ is zero by Remark~\ref{degree-two}.
We consider three cases, depending on whether the spaces
$(V^g)^\perp$ and $(V^h)^\perp$ intersect in dimension $0$,\, $1$, or $2$.

{\bf Case 1: Disjoint orthogonal complements.}
 Assume
that $(V^g)^\perp\cap (V^h)^\perp =0$.  Then
$$\codim V^g + \codim V^h = \codim V^{gh}=\codim V^{hg}$$
(e.g., see~\cite[Lemma~2.1]{paper1}). 
Thus $\codim V^{gh}=4 = \codim V^{hg}$, and by examination of 
Definition~\ref{H,Hg},
$H^3_{gh}=0 = H^3_{hg}$ (as $H^{\DOT}_{gh}$ ``begins'' in degree 4).
But we have excluded this case from consideration.

{\bf Case 2. Equal orthogonal complements.}
Assume that $(V^g)^\perp$ and $(V^h)^\perp$ are equal, i.e.,
$\dim ((V^g)^{\perp}\cap (V^h)^{\perp})=2$, and $v_1, v_2$ span $(V^g)^\perp=(V^h)^\perp$.
Consider
$$\alpha=f_g \bar{g} \otimes (v_1^*\wedge v_2^*)
\quad\text{and}\quad
\beta =f_h' \bar{h} \otimes (v_1^*\wedge v_2^*),$$ where
$f_g, f_h'\in S(V^g)=S(V^h)$.
If nonzero,  
$
\ttau_1{\alpha}\bigl(\ttau_2{\beta}
(v_{i}\otimes v_{j}) \otimes {}^hv_{k}\bigr) 
=
\ttau_1{\alpha}\bigl(f_h' \otimes\ ^hv_{k}\bigr)
$
up to a constant.
But $\ttau_1{\alpha}\bigl(f_h' \otimes *\bigr)=0$
by Remark~\ref{zeroterm}, as $f_h'\in S(V^g)$ has zero partial derivative
with respect to any basis element in $(V^g)^\perp$.
Hence 
$$0=\ttau_1{\alpha}\bigl(\ttau_2{\beta}
(v_{i}\otimes v_{j}) \otimes {}^hv_{k}\bigr).
$$
Similarly,
$$0\, =\, 
\ttau_1{\alpha}\bigl(v_{i} \otimes
\ttau_2{\beta}(v_{j}\otimes v_{k})\bigr)
\,  = \,
\ttau_2{\beta}\bigl(\ttau_1{\alpha}(v_{i} \otimes 
v_{j})\otimes {}^gv_{k}\bigr) 
\, = \, 
\ttau_2{\beta}\bigl(v_{i}\otimes\ttau_1{\alpha}
(v_{j}\otimes v_{k})\bigr)
$$
for all distinct $i$, $j$, and $k$ and thus
Equation~\ref{four2} yields zero.  Hence, $[[\alpha,\beta]]=0$.

{\bf Case 3. Overlapping orthogonal complements.}
Assume the two spaces 
$(V^g)^{\perp}$ and $(V^h)^{\perp}$ overlap only partially, i.e., 
$\dim ((V^g)^{\perp}\cap (V^h)^{\perp})=1$.
The remainder of the proof is devoted to this last case.

We now refine our bases $B_1$ and $B_2$
to ease determination of the prebracket.  
Recall that
$B_1=\{w_1,\ldots, w_n\}$
and $B_2=\{v_1,\ldots, v_n\}$
are bases of $V$ with $w_1,w_2$
in $(V^g)^\perp$, $v_1,v_2$ in $(V^h)^{\perp}$, 
$w_3,\ldots, w_n$ in $V^g$, and $v_3,\ldots, v_n$ in $V^h$.
We make additional assumptions.
Let $W=(V^g)^\perp+(V^h)^\perp$.  The space $W$ has dimension $3$
since $(V^g)^\perp$ and $(V^h)^\perp$ intersect
in dimension 1.  Notice that $V^g$ intersects
$W$ nontrivially, 
as otherwise $\dim(W+V^g)=\dim(W)+\dim(V^g)=3+(n-2)=n+1>\dim V$.
Also note that $W^\perp =V^g \cap V^h$.  
Thus, we may assume further that
$w_1, w_2, w_3$ span $W$, $w_4, \ldots, w_n$ span $W^{\perp}$,
and $v_4=w_4, \ldots, v_n=w_n$.

These refining assumptions have no effect on the prebracket $[[\alpha,\beta]]$.
Although $\ta$ is independent of choice of bases as a map on cohomology,
it depends on the choices $B_1$ and $B_2$ as a {\em cochain map}.
(See cautionary Remark~\ref{caution}.)
Yet in refining our choices of $B_1$ and $B_2$, we have {\em not} altered the
values of $\ta(\alpha)$ or $\ta(\beta)$ as cochains.
Indeed, quantum partial differentiation with respect to one subset of variables
in a basis ignores any change of basis affecting only the other variables.
The exterior part of $\alpha$ is an element of $\bigwedge ((V^g)^\perp)^*$
(as $\alpha$ lies in $H_g^2$),
so the map $\ta(\alpha)$ differentiates with respect to 
vectors in $(V^g)^\perp$.  Since we altered the basis of $B_1$ on
$V^g$ alone, the map
$\ta(\alpha)$ is unchanged. Similarly, $\ta(\beta)$ is also
unchanged.

We examine the coefficient of $\overline{gh}$ in Equation~\ref{four2}.
Consider
$$\alpha=f_g \bar{g} \otimes w_1^*\wedge w_2^*\quad\text{and}\quad
\beta=f_h' \bar{h} \otimes v_1^*\wedge v_2^*,$$ where
$f_g \in S(V^g)$ and
$f_h' \in S(V^h)$.
We rewrite the sum (giving the coefficient of ${\overline{gh}}$)
with indices in a different order:
\begin{equation}
\label{simplifiedfour}
\begin{aligned}
\sum_{\pi\in \Sym_3}\sgn(\pi)\, 
&
\biggl[
\ttau_1{\alpha}\Bigl(\ttau_2{\beta}
(v_{\pi(i)}\otimes v_{\pi(j)}) \otimes\, ^h v_{\pi(k)}\Bigr) 
-  
\ttau_1{\alpha}\Bigl(v_{\pi(i)} \otimes
\ttau_2{\beta}(v_{\pi(j)}\otimes v_{\pi(k)})\Bigr)
\biggr]\\
=
\sum_{\pi\in \Sym_3}\sgn(\pi)\, 
&
\biggl[
\ttau_1{\alpha}\Bigl(\ttau_2{\beta}
(v_{\pi(i)}\otimes v_{\pi(j)}) \otimes\, ^h v_{\pi(k)}\Bigr) 
- 
\ttau_1{\alpha}\Bigl(v_{\pi(k)} \otimes
\ttau_2{\beta}(v_{\pi(i)}\otimes v_{\pi(j)})\Bigr)
\biggr]\ .
\end{aligned}
\end{equation}
By Remark~\ref{zeroterm}, each of
the above summands is zero for every permutation save one.
Indeed, the summand corresponding to $\pi$ is zero unless
$\pi(i)=1$, $\pi(j)=2$, and $\pi(k)=3$
(since $v_4=w_4,\ldots, v_n=w_n$), and we are left with
\begin{equation}
\label{simplifiedtwo}
\begin{aligned}
 &\ttau_1{\alpha}\Bigl(\ttau_2{\beta}
(v_{1}\otimes v_{2}) \otimes\, ^hv_{3}\Bigr) 
\ - \ 
\ttau_1{\alpha}\Bigl(v_{3} \otimes
\ttau_2{\beta}(v_{1}\otimes v_{2})\Bigr)\\
= \
& \ttau_1{\alpha}\Bigl(f_h' \otimes v_{3}\Bigr) \ - \ 
\ttau_1{\alpha}\Bigl(v_{3} \otimes f_h'\Bigr) \ .
\end{aligned}
\end{equation}

We show that this difference maps to zero under the projection
$S(V)\rightarrow S(V^{gh})$.  
We reduce to the case when $f_h'$ is a power of $v_3$.
Indeed, the above difference is just
$$
\del_1(f_h')\ ^{s_1} \del_2 (v_{3}) -  
\del_1(v_{3})\ ^{s_1} \del_2(f_h')
$$ 
where $s_1$ is some reflection
and $\del_i$ is some partial quantum differentiation with respect to $w_i$
in the basis $\{w_1, w_2, w_3, w_4\text{$=$}v_4, \ldots, w_n\text{$=$}v_n\}$.
As $\del_1,\del_2$ are both $\CC[v_4,\ldots, v_n]$-linear, 
and $f_h'$ lies in $S(V^h)=\CC[v_3, \ldots, v_n]$,
we may break $f_h'$ into
its monomial summands (in the basis $\{v_1,\ldots, v_n\}$)
and pull out all factors from $\CC[v_4, \ldots, v_n]$
when evaluating the above difference.
Thus, it suffices to consider
the special case when $f_h'=v_3^m$ for some $m\geq 0$.

By Lemma~\ref{udivides}, $u:=v_3- {}^gv_3$ divides
\begin{equation}
\label{difference}
\begin{aligned}
\ttau_1{\alpha}(v_3^m \otimes v_{3}) \ - \ 
\ttau_1{\alpha}(v_{3} \otimes v_3^m ) \ . 
\end{aligned}
\end{equation}
But notice that $u$ lies in $(V^{gh})^{\perp}$, since
$$
  ^{gh}(u) = \ ^{gh}(v_3-\, ^{g}v_3) =\ ^{g}v_3-\, ^{ghg}v_3=\ ^{1-gh}(\, ^gv_3),
$$
i.e., $^{gh}u$ lies in ${\rm im}(1-gh)=(V^{gh})^{\perp}$.
Thus the difference~(\ref{difference}) 
lies in the ideal $I((V^{gh})^{\perp})$ of $S(V)$
and projects to zero under the map $S(V)\rightarrow S(V^{gh})$.

By a symmetric argument, the coefficient of $\overline{hg}$
in Equation~\ref{four2} projects to zero 
under the map $S(V)\rightarrow S(V^{hg})$.
Thus, the prebracket $[[\alpha, \beta]]$
projects to zero under the map $\projH$.
The same arguments apply to 
$[[\,^a\alpha,\,^b\beta\,]]_{\rule[0ex]{0ex}{2ex}(\,^aB_g,\, ^bB_h)}$ 
for arbitrary $a,b$ in $G$.
Hence,
$$
[\alpha, \beta]
:= \frac{1}{|G|^2}
\ \projH
\sum_{\substack{a, b \in G }}
[[\, ^{a}\alpha,\, ^{b}\beta]]_{\rule[0ex]{0ex}{2ex}(\,^aB_g,\, ^bB_h)}
\ = \ 0\ .
$$
\end{proof}
\vspace{1ex}
%%%%%%%%%%%%%%%%%%%%%%%%%%%%%%%%%%%%%%%%%%%%%%%%%%%%%%%%%%%%%%%%%%%%%%

The theorem above implies that if $\alpha$ lies in $\HH^2(S(V)\# G)$ with
$[\alpha,\alpha]\neq 0$, then the support of $\alpha$ includes at
least one group element acting as the identity on $V$:
 
%%%%%%%%%%%%%%%%%%%%%%%%%%%%%%%%%%%%%%%%%
%%%%%%%%%%%%%%%%%%%%%%%%%%%%%%%%%%%%%%
\vspace{2ex}
\begin{cor}
\label{poissonstructures}
The Gerstenhaber square bracket of every $\alpha$ in $\HH^2(S(V)\#G)$ supported
off of $K$ is zero, i.e., $\alpha$ defines a noncommutative Poisson structure
on $S(V)\# G$:
$$
[\alpha,\alpha]=0\ .
$$ 
\end{cor}
\vspace{2ex}
%%%%%%%%%%%%%%%%%%%%%%%%%%%%%%%%%%%%%

Next, we illustrate our results by 
giving an explicit example of a Gerstenhaber bracket in degree 2.

\vspace{3ex}
\begin{example}\label{examples}
Let $G=D_8$, the dihedral group of order 8, generated by $g$ and $h$
with relations $g^4=1=h^2$, $hgh^{-1}=g^3$, 
realized as a subgroup of  ${\rm {GL}}_3(\CC)$ in the
following way:
$$
  g = \left(\begin{array}{ccc}0&i&0\\
                              i&0&0\\
                              0&0&1
            \end{array}\right) , \ \ \ 
  h = \left(\begin{array}{ccc} 1 & 0 & 0\\
                               0 & -1& 0\\
                               0 & 0 & -1
             \end{array}\right),
$$
where $i=\sqrt{-1}$.
Let $B_h=\{v_1,v_2,v_3\}$ be the corresponding basis of $\CC^3$,
and set $B_g=\{w_1,w_2, w_3\}$ where
$w_1=v_1+v_2$, $w_2=-v_1+v_2$, and $w_3=v_3$, so that 
${}^gw_1=iw_1$ and ${}^gw_2=-iw_2$.
Define 
\begin{eqnarray*}
  \alpha &=& v_3 \, \overline{g}\ot w_1^*\wedge w_2^* 
\quad\quad\text{in}\ H^2_g\\
  \beta &=& v_1^3\, \overline{h}\ot v_2^*\wedge v_3^* 
\quad\quad\ \text{in}\ H^2_h\ .
\end{eqnarray*}
Note that $v_3\in S(V^g)$, $v_1^3\in S(V^h)$, and $\alpha$,
$\beta$ are $Z(g)$-, $Z(h)$-invariant, respectively.
We compute $[\alpha,\beta\hspace{.1ex}]$.
By Definition~\ref{prebracketformula}, 
$(\alpha\overline{\circ}\, \beta)(v_{1}\wedge v_{2}\wedge v_{3})$ is equal to
$$
    \sum_{\pi\in\Sym_3}(\sgn \pi)\ttau\alpha \bigl( 
    \ttau{\beta}(v_{\pi(1)}\ot v_{\pi(2)})\ot {}^hv_{\pi(3)}\bigr)
              \overline{gh}
   -\ttau\alpha\bigl(v_{\pi(1)}\ot \ttau{\beta}(v_{\pi(2)}\ot
    v_{\pi(3)})\bigr)\overline{gh}.
$$
The $\pi = 1$ and $\pi = (123)$ terms are 
\begin{eqnarray*}
&&
\biggl(-\frac{3}{16}w_1^2 -\frac{3}{16}(1+i)w_1w_2-\frac{i}{16}w_2^2\biggr)\,
v_3\ \overline{gh}
       \ \mbox{ and}\\
&&
\biggl(-\frac{i}{16}w_1^2+\frac{3}{16}(1+i)w_1w_2-\frac{3}{16}w_2^2\biggr)\, 
v_3 \ \overline{gh},
\end{eqnarray*}
respectively, and the remaining terms are $0$.
Combining these, we find that the
coefficient of $\overline{gh}$ in 
$(\alpha\overline{\circ}\beta)(v_1\wedge v_2\wedge v_3)
$ is 
$$
\left(\frac{-3-i}{16}w_1^2+\frac{-3-i}{16}w_2^2\right)v_3.
$$
Thus $\alpha\overline{\circ}\beta$ is nonzero, as a function at the chain level,
given our choices of bases. By Lemma~\ref{udivides},
$(\alpha\overline{\circ}\beta )(v_1\wedge v_2\wedge v_3)$ is divisible by
$u:=(1-g)v_1 = v_1-iv_2$, and indeed 
$w_1^2+w_2^2$ is a scalar multiple of $(v_1+iv_2)(v_1-iv_2)$, the
product of an element in $V^{gh}$ with an element in $(V^{gh})^{\perp}$.
Hence, $[\alpha,\beta]$ projects to zero under the map $\projH$, i.e.,
it is a coboundary.
Similar calculations give 
${}^{a}\alpha\,\overline{\circ}\ ^{b}\beta$ and ${}^b\beta\,\overline{\circ}\
{}^a\alpha$ for all $a,b$ in $G$.
Thus $[\alpha,\beta]=0$ in cohomology, as predicted by Theorem~\ref{bracketalwayszero}.
\end{example}
\vspace{3ex}

%%%%%%%%%%%%%%%%%%%%%%%%%%%%%%%%%%%%%%%%%%%%%%%%%%%%%%%%%%%%%%%%%%%%%%
%%%%%%%%%%%%%%%%%%%%%%%%%%%%%%%%%%%%%%%%%%%%%%%%%%%%%%%%%%%%%%%%%%%%%%
%%%%%%%%%%%%%%%%%%%%%%%%%%%%%%%%%%%%%%%%%%%%%%%%%%%%%%%%%%%%%%%%%%%%%%
\section{Abelian groups and inner products of characters}
\label{sec:abelian}

We now consider the Gerstenhaber bracket for abelian groups.
Orthogonality relations on characters allow us to place
Theorem~\ref{bracketalwayszero} in context.  We observe
that the hypothesis of the theorem cannot be weakened and 
we show that its converse is false.

Let $G$ be a finite abelian group and $V$ a 
(not necessarily faithful) $\CC G$-module 
of finite dimension $n$.
We first explain how inner products of characters of $G$
determine the Gerstenhaber bracket on $\HHD(S(V)\# G)$.
It would be interesting to know whether a similar description holds
for arbitrary groups.  
We concentrate on cohomological degree $2$ due to  
connections with deformation theory.
Assume $\dim V \geq 3$, as otherwise the Hochschild cohomology of 
$S(V)\# G$ in degree 3 is always zero.

Fix a simultaneous basis of eigenvectors for $G$,
say $B=\{v_1,\ldots,v_n\}$.
For $i\in\{1,\ldots,n\}$, let $\chi_i$ be the character
defined by ${}^{g}v_i = \chi_i(g)\, v_i$  for each $g$ in $G$.
The Gerstenhaber bracket of $2$-cocycles in 
$(H^{\DOT})^G\cong\HHD(S(V)\#G)$
may be computed summand by summand, at the cochain level,
as the bracket is linear and extends to $H^{\DOT}$
(see Remark~\ref{extend}). 
Thus, it suffices to describe the bracket
on simple cochains of the form 
$$\alpha=f_g \overline{g} \ot v_i^*\wedge v_j^*
\quad\text{and}\quad 
\beta=f_h' \overline{h} \ot v_k^*\wedge v_m^*$$
where $f_g, f_h'$ are monomials and $g,h$ lie in $G$.  Three cases arise.

{\bf Case 1:}
If the exterior parts of $\alpha$ and $\beta$
agree up to sign (i.e., $\{i,j\}=\{k,m\}$), then the bracket
$[\alpha,\beta]$ (given by Equation~(\ref{special}))
is easily seen to be zero
by Remark~\ref{zeroterm}. 

{\bf Case 2:}
If the exterior parts of $\alpha$ and $\beta$
partially overlap, relabel indices so that
$$\begin{aligned}  
\alpha =
f_g\,\overline{g}\ot v_1^*\wedge v_2^* 
\quad\quad\text{and}\quad\quad
\beta = 
f_h'\,\overline{h}\ot v_2^*\wedge v_3^*\ .
\end{aligned}
$$
The bracket operations $[*, \beta]$ and $[\alpha, *]$ are then
both $\CC[v_4,\ldots, v_n]$-linear: 
$[v_i\alpha,\beta]=[\alpha,v_i\beta]=v_i[\alpha,\beta]$ for $i\geq 4$.
The bracket computation thus reduces to that of Example~\ref{abelian},
in which $f_g=v_1^{c_1}v_2^{c_2}v_3^{c_3}$ and 
$f_h'=v_1^{d_1}v_2^{d_2}v_3^{d_3}$.
The formula in that example
expresses the bracket $[\alpha,\beta]$ 
in terms of inner products of the group characters $\chi_i$.
Note that the bracket $[\alpha,\beta]$ is zero or
$\chi_1^{c_1+d_1-1}\chi_2^{c_2+d_2-2}\chi_3^{c_3+d_3-1}=1$,
and thus the bracket $[\alpha,\beta]$ is $G$-invariant (see the formula),
even if $\alpha$ and $\beta$ are not
(verifying Remark~\ref{extend}).

{\bf Case 3:}
If the exterior parts of $\alpha$ and $\beta$ do not overlap
(i.e., $\{i,j\}\cap\{k,m\}=\emptyset$),
a similar but lengthier formula results.
We may assume
$$\begin{aligned}  
\alpha =
(v_1^{c_1}v_2^{c_2}v_3^{c_3}v_4^{c_4})\,
\overline{g}\ot v_1^*\wedge v_2^* 
\quad\text{and}\quad
\beta = 
(v_1^{d_1}v_2^{d_2}v_3^{d_3}v_4^{d_4})
\,\overline{h}\ot v_3^*\wedge v_4^*\ .
\end{aligned}
$$
The bracket $[\alpha,\beta]$ is given by
$$
\Proj_H 
\sum_{1\leq l\leq 4} \kappa_l\,
\bigl(v_1^{c_1+d_1}v_2^{c_2+d_2}v_3^{c_3+d_3}v_4^{c_4+d_4}v_l^{-1}\bigr)
\, \overline{gh}\otimes 
v_1^*\wedge\cdots\wedge\widehat{v_l^*}\wedge\cdots\wedge v_4^*\ ,
$$
where 
$$
\begin{aligned}
&\kappa_1:= 
+
\langle
\chi_1^{c_1-1}\chi_2^{c_2+1}\chi_3^{c_3}\chi_4^{c_4}
, 1 \rangle
\langle
\chi_1^{d_1}\chi_2^{d_2-2}\chi_3^{d_3-1}\chi_4^{d_4-1}
, 1 \rangle
\cdot
[d_1]_{\chi_1(g)},
\\
&\kappa_2:= 
-
\langle
\chi_1^{c_1+1}\chi_2^{c_2-1}\chi_3^{c_3}\chi_4^{c_4}
, 1 \rangle
\langle
\chi_1^{d_1-2}\chi_2^{d_2}\chi_3^{d_3-1}\chi_4^{d_4-1}
, 1 \rangle
\cdot
[d_2]_{\chi_2(g)}
\cdot\chi_1^{d_1}(g),
\\
&\kappa_3:= 
+
 \langle
\chi_1^{c_1-1}\chi_2^{c_2-1}\chi_3^{c_3}\chi_4^{c_4-2}
, 1 \rangle
\langle
\chi_1^{d_1}\chi_2^{d_2}\chi_3^{d_3-1}\chi_4^{d_4+1}
, 1 \rangle
\cdot[c_3]_{\chi_3(h)}
\cdot(\chi_1^{c_1}\chi_2^{c_2})(h),
\\
& \kappa_4:=
-
\langle
\chi_1^{c_1-1}\chi_2^{c_2-1}\chi_3^{c_3-2}\chi_4^{c_4}
, 1 \rangle
\langle
\chi_1^{d_1}\chi_2^{d_2}\chi_3^{d_3+1}\chi_4^{d_4-1}
, 1 \rangle 
\cdot[c_4]_{\chi_4(h)}\cdot
(\chi_1^{c_1}\chi_2^{c_2}\chi_3^{c_3})(h),
\end{aligned}
$$
$\langle\ , \ \rangle$ denotes inner product of 
group characters and
$[m]_\lambda$ is the quantum integer
$1+\lambda+\lambda^2+\ldots+\lambda^{m-1}$
(or zero when $m=0$).
%each constant $\kappa_m$ is determined by 
%inner products of group characters.
%We do not record the precise formula here, but merely 
Note that
$[\alpha,\beta]$ is either zero or
$\chi_1^{c_1+d_1-1}\chi_2^{c_2+d_2-1}\chi_3^{c_3+d_3-1}\chi_4^{c_4+d_4-1}=1$.
Thus the bracket $[\alpha,\beta]$ is $G$-invariant,
even if $\alpha$ and $\beta$ are not, as results in previous sections predict.

\medskip

We now use our analysis of the Gerstenhaber bracket for abelian groups 
to revisit Theorem~\ref{bracketalwayszero} and Corollary~\ref{poissonstructures}.
The next two results show that the hypotheses on 
Theorem~\ref{bracketalwayszero} and Corollary~\ref{poissonstructures}
cannot be weakened.  Theorem~\ref{nonzeroabeliansquarebracket} below reveals that 
the converse of Theorem~\ref{bracketalwayszero} is false for any abelian group.
We do not know whether the statements below hold for nonabelian groups.
Recall that $K=\{g\in G: V^g=V\}$, the kernel of the action of $G$ on $V$.

%%%%%%%%%%%%%%%%%%%%%%%%%%%%%%%%%%%%%%%%%%%%%%%%%%%%
\vspace{2ex}
\begin{prop}\label{nonzeroabelianbracket}
Let $G$ be an abelian group.  Suppose $g,h$ lie in $K$ and $\dim V $
is at least $3$.
Then there are elements $\alpha, \beta$ in $\HH^2(S(V)\#G)$
supported on $g,h$, respectively,
with nonzero Gerstenhaber bracket:
$[\alpha,\beta]\neq 0$.
\end{prop}
\vspace{1ex}
\begin{proof}
We shall use the notation and formula of Example~\ref{abelian}.
Let 
$$
\begin{aligned}
\alpha = \bigl(v_1v_2^{|G|+1} \bigr)\, \overline{g} \otimes v_1^*\wedge v_2^*\ \quad\text{and}\quad
\beta = \bigl(v_2 v_3^{|G|+1}\bigr)\, \overline{h}\ot v_2^*\wedge v_3^*\ ,
\end{aligned}
$$ 
i.e., set
$c_1=1$, $c_2= |G|+1$, $c_3 =0$, $d_1=0$, $d_2=1$, and $d_3= |G|+1$.
Note that $\alpha$ and $\beta$ are invariant cocycles.
Then
$$  
  [\alpha,\beta] =
   \kappa\   \bigl(v_1 v_2^{|G|+1} v_3^{|G|+1}\bigr)\, \overline{gh}
         \otimes v_1^*\wedge v_2^*\wedge v_3^* 
$$
where
$$
\kappa =
   \langle 1\cdot 1 , \ 1 \rangle \langle 1, \ 1\cdot 1\rangle 
   (|G|+1 -1) \, 
     = \ |G| \ 
 $$
and $1$ denotes the trivial character of $G$.
Hence $[\alpha,\beta]$ is nonzero.
\end{proof}

%%%%%%%%%%%%%%%%%%%%%%%%%%%%%%%%%%%%%%%%%%%%%%%%%%

The last proposition implies:

%%%%%%%%%%%%%%%%%%%%%%%%%%5
\vspace{2ex}
\begin{thm}\label{nonzeroabeliansquarebracket}
Let $G$ be an abelian group. Then
\begin{itemize}
\item
There is a 2-cocycle supported on $K$ whose square bracket is nonzero.
\item
There is a 2-cocycle supported on $K$ whose square bracket is zero.
\end{itemize}
\end{thm}
\vspace{1ex}
\begin{proof}
We prove a slightly stronger statement.  Let $k$ be any element of $K$.
We apply Proposition~\ref{nonzeroabelianbracket} in the case that $g=h=k$.
We obtain cocycles $\alpha, \beta$ in $(H_k^{\DOT})^G\cong\HH^2(S(V)\#G)$ with
$[\alpha,\beta]\neq 0$. Then as
$$[\alpha + \beta, \alpha+\beta]=[\alpha, \alpha]
+2[\alpha,\beta]+[\beta,\beta]\ ,$$
there must be a cocycle supported on $k$
with square bracket nonzero.

Now set $\alpha=(v_1 v_2)\overline{k}\otimes v_1^*\wedge v_2^*$.
Since $G$ acts diagonally, $\alpha$ is $G$-invariant, i.e., $\alpha$ lies
in $(H^{2})^G\cong\HH^2(S(V)\#G)$. 
Yet $[\alpha,\alpha]=0$ (see Case 1 above).
\end{proof}
%%%%%%%%%%%%%%%%%%%%%%%%%%%%%%5
\vspace{1ex}

%%%%%%%%%%%%%%%%%%%%%%%%%%%%%%

We end this section by pointing out 
a direct and easy proof of 
Theorem~\ref{bracketalwayszero} for abelian groups as follows.
Suppose $\alpha$ and $\beta$ in $(H^2)^G$ are supported off $K$ but
$[\alpha, \beta]\neq 0$.
Then the bracket of some summand of $\alpha$ and some summand of $\beta$
is nonzero.  
We consider the three cases at the beginning of this section.
The bracket in Case 1 is always zero.
The bracket in Case 3 is also zero:
Remark~\ref{degree-two} implies that 
$v_1,v_2$ span $(V^g)^{\perp}$
while $v_3, v_4$ span $(V^h)^{\perp}$; hence
$\CC[v_1,v_2,v_3,v_4]$ projects to zero under the map
$S(V)\rightarrow S(V^{gh})=S(V^{hg})$ and the bracket lies
in the kernel of $\Proj_H$.
We thus reduce to Case~2 and Example~\ref{abelian}: 
We assume $v_1,v_2$ span
$(V^g)^{\perp}$ and $v_2, v_3$ span 
$(V^h)^{\perp}$.
Then $[\alpha,\beta]=0$ as the polynomial coefficient of $\alpha$
(in $S(V^g)$) and the polynomial coeffient of $\beta$ (in $S(V^h)$)
both vanish
after taking the partial derivative with respect to $v_2$.
(See Remark~\ref{zeroterm} or 
the formula of Example~\ref{abelian} with $c_1=c_2=d_2=d_3=0$.)

%%%%%%%%%%%%%%%%%%%%%%%%%%%%%%%%%%%%%%%%%%%%%%%%%%%%%%%%%%%%%%
%%%%%%%%%%%%%%%%%%%%%%%%%%%% SECTION %%%%%%%%%%%%%%%%%%%%%%%%%
%%%%%%%%%%%%%%%%%%%%%%%%%%%%%%%%%%%%%%%%%%%%%%%%%%%%%%%%%%%%%%

%%%%%%%%%%%%%%%%%%%%%%%%%%%%%%%%%%%%%%%%%%%%%%%%%%%%%%%%%%%%%%%%%%%%%%555
%%%%%%%%%%%%%%%%%%%%%%%%%%%%%%%%%%%%%%%%%%%%%%%%%%%%%%%%%%%%%%%%%%%%%%%%%5

%%%%%%%%%%%%%%%%%%%%%%%%%%%%%%%%%%%%%%%%%%%%%%%%%%%%%%%%%%%%
%%%%%%%%%%%%%%%%%%%%%%%%%%%%%%%%%%%%%%%%%%%%%%%%%%%%%%%%%%%%
\section{Graded (Drinfeld) Hecke algebras}\label{previouspaper}
We end by briefly highlighting connections with graded Hecke algebras
(or Drinfeld Hecke algebras),
which include symplectic reflection algebras
(and rational Cherednik algebras).
We show how the maps in previous sections give explicit conversions
among graded Hecke algebras, deformation theory,
and Hochschild cocycles (expressed as vector forms).

Let $G$ be a finite group and let $V$ be a finite dimensional $\CC G$-module. 
Let $T(V)$ denote the tensor algebra on $V$
and let $\kappa:V\times V\rightarrow \CC G$
be a bilinear, skew-symmetric function.
A {\bf graded Hecke algebra} is a quotient 
\begin{equation}\label{gHa}
   \mathcal{H}=
   T(V)\# G / (v\ot w-w\ot v-\kappa(v,w)\,\mid v,w\in V)\ 
\end{equation}
that satisfies the Poincar\'e-Birkhoff-Witt property:
Any linear splitting of the canonical projection $T(V)\rightarrow S(V)$
induces a vector space isomorphism $\mathcal{H}\cong S(V)\# G$.
We extend scalars to the polynomial ring $\CC[t]$ and consider 
every graded Hecke algebra as a quotient
\begin{equation}\label{hecke-relations}
  T(V)\# G [t]/(v\ot w-w\ot v - \kappa(v,w)\,t\ |\ v,w\in V)\ .
\end{equation}
See~\cite{RamShepler,SheplerWitherspoon} for basic definitions.

We first restate~\cite[Theorem~3.2]{Witherspoon2} in our context.
(Alternatively, results from~\cite{BGS} and~\cite{BravermanGaitsgory}
could be used to obtain this and related statements on 
deformations of $S(V)\# G$.)
Recall that the $i$-th multiplication map for a given deformation of $S(V)\# G$
is denoted $\mu_i$ (see Section~\ref{Prelim}). 
Consider $S(V)\# G$ to be a graded algebra with $\deg v=1$, $\deg g = 0$ for
all $v$ in $V$, $g$ in $G$.  We agree that the zero map 
has degree $i$ for any integer $i$.

%%%%%%%%%%%%%%%%%%%%%%%%%%%%%%%%%%%%%%%%%%%%%%5
\vspace{2ex}
\begin{thm}\label{degree-mu-i}
Every graded Hecke algebra is isomorphic to a deformation of $S(V)\# G$.
In fact, up to isomorphism, the graded Hecke algebras are precisely the deformations
of $S(V)\# G$ over $\CC[t]$
for which the $i$-th multiplication map
lowers degree by $2i$, i.e., $\mu_i$ is a (homogeneous) graded map with 
$\deg \mu_i = -2i$ ($i\geq 1$).
\end{thm}
\vspace{2ex}
%%%%%%%%%%%%%%%%%%%%%%%%%%%%%%%%%%%%%%%%%%%%%%%%%%%%%%%%%%%%%%

We now discuss the explicit conversions among graded Hecke algebras,
deformations, and Hochschild cocycles by interpreting the above
theorem and its proof using our results from previous sections.

%%%%%%%%%%%%%%%%%%%%%%%%%%%%%%%%%%%%%%%%%%%%%%%%%%%5
\subsection*{Deformations to Graded Hecke Algebras}
Given a fixed deformation of $S(V)\# G$ over $\CC[t]$ for which
$\deg\mu_i = -2i$ ($i\geq 1$), we obtain a graded Hecke algebra 
by defining $\kappa: V\times V \rightarrow \CC G$ by
$$
   \kappa(v,w) = \mu_1(v\ot w) - \mu_1(w\ot v).
$$
Note that $\kappa$ is skew-symmetric by definition, even if
$\mu_1$ is not skew-symmetric. It is shown in~\cite{Witherspoon2}
that the graded Hecke algebra (\ref{gHa}) corresponding to this choice
is isomorphic to the deformation with which we started.

%%%%%%%%%%%%%%%%%%%%%%%%%%%%%%%%%%%%%%%%%%%%%%
\subsection*{
Hochschild Cocycles to Graded Hecke Algebras}

We identify 
$\HH^{\DOT}(S(V)\#G)$ with a subset of $C^{\DOT}$,
vector forms tagged by group elements, using Theorem~\ref{Decomposition}.
The algebra $S(V)$ is graded by polynomial degree:
$S(V)=\bigoplus_{k\geq 0} S(V)_{(k)}$.
This induces a grading on $S(V)\# G$ (after assigning degree 0 to each
$g$ in $G$) which is inherited by $C^{\DOT}$:
$$C^{\DOT}=\bigoplus_{k\geq 0,\ g\in G}
S(V)_{(k)}\, \overline{g}\otimes \Wedge^{\DOT}V^*.$$
Recall that a Hochschild $p$-cocycle is said to be {\em constant} if it
lies in the $0$-th graded piece of $C^{p}$, i.e.,
defines a vector form with constant polynomial part.
We rephrase Theorem~8.7 of~\cite{SheplerWitherspoon},
which uses Theorem~\ref{degree-mu-i} to determine
that every graded Hecke algebra
arises from a {\em constant} $2$-form.

%%%%%%%%%%%%%%%%%%%%%%%%%%%%%
\vspace{2ex}
\begin{theorem}\label{constantcocycles}
The parameter space of graded Hecke algebras is isomorphic to
the space of {\em constant} Hochschild $2$-cocycles,
$$
\biggl(\,\bigoplus_{g\in G} 
    \overline{g}\ot \Wedge^{2 - \codim V^g}(V^g)^*
     \ot \Wedge^{\codim V^g}((V^g)^{\perp})^*\biggr)^{G}.
$$
\end{theorem}
\vspace{2ex}
%%%%%%%%%%%%5

We next give an explicit conversion from constant 2-cocycles to
graded Hecke algebras. 

%%%%%%%%%%%%%%%%%%%%%%%%%%%%%%%%%%%%%%%%%%5
\vspace{2ex}
\begin{prop}
The correspondence above
is induced from the  map:
$$\begin{aligned}
\{ \text{constant Hochschild $2$-cocycles} \} 
&\rightarrow 
\{\text{graded Hecke algebras}\}\\
\alpha\ \ \ &\mapsto\ \ \  (T(V)\#G)/I_\alpha\ ,
\end{aligned}$$ 
where $I_\alpha$ is the ideal
generated by
$\{v \ot w-w \ot v - \alpha(v\wedge\,w)| \ v, w \in V\}$
and $\alpha \in C^2\cong\Hom_{\CC}(\Wedge^2 V, \CC G)$ .
(I.e., $\alpha$ defines a graded Hecke
algebra with $\kappa(v,w)=\alpha(v\wedge w)$.)
\end{prop}
%%%%%%%%%%%%%%%%%%%%%%%%%%%%%%%%%%%%%%%%%
\vspace{1ex}
\begin{proof}
We use our conversion map
$\gamm:(H^{\DOT})^G\rightarrow \HHD(S(V)\#G)$ of Theorem~\ref{S:CGW2}.
The values of the multiplication map
$\mu_1$ in a deformation are given by the corresponding Hochschild 2-cocycle
in $\HHD(S(V)\#G)$.  Let $\alpha$ lie in $(H^2)^G$ with isomorphic
image $\gamm(\alpha)$ in $\HHD(S(V)\#G)$.
Note that for any $v,w$ in $V$,
$$(\ta\alpha)(v\ot w - w\ot v)=\alpha(v\wedge w)\ .$$
Then,
$$
\begin{aligned}
\gamm(\alpha)(v\ot w-w\ot v)
& =\frac{1}{|G|} \sum_{g\in G} \rh^*\ ^g(\ta\alpha)(v\ot w-w\ot v)\\
& =\frac{1}{|G|} \sum_{g\in G}\ ^g\bigl((\ta\alpha)\ ^{g^{-1}}(v\ot
w-w\ot v)\bigr)\\
& =\frac{1}{|G|} \sum_{g\in G}\ ^g\bigl(\alpha\ ^{g^{-1}}(v\wedge
w)\bigr)\\
& =\frac{1}{|G|} \sum_{g\in G}\ (\, ^g\alpha)(v\wedge w)\\
& = \alpha(v\wedge w)\ .
\end{aligned}
$$
By Theorem~\ref{degree-mu-i} and its proof,
$\mu_i(v\ot w)=0$ for all $i\geq 2$.
Thus we have 
$\gamm(\alpha)(v\ot w - w\ot v)t = (\mu_1(v\ot w) - \mu_1(w\ot v))t
= v* w -w*v$ 
in the deformation of $S(V)\#G$ over $\CC[t]$. This
corresponds to the relation in the graded Hecke algebra
$\alpha(v\wedge w)= v\ot w - w \ot v$.  
\end{proof}
\vspace{1ex}

%%%%%%%%%%%%%%%%%%%%%%%%%%%5
\subsection*{Hochschild Cocycles to
Infinitesimal Deformations}
The interpretation in the last result and its proof  hold for 
Hochschild cocycles
of arbitrary polynomial degree, not just constant cocycles:
The proof gives the conversion from
{\em any} Hochschild $2$-cocycle to an infinitesimal deformation.
In fact, a closed form expresses the multiplication map $\mu_1$ 
in terms of quantum differentiation.
Let $\alpha$ in $(H^{\DOT})^G$ be a Hochschild $2$-cocycle;
then $\alpha$ defines a multiplication map $A\otimes A\rightarrow A$
for $A=S(V)\#G$ given by
\begin{equation}\label{productoftwo}
\mu_1(f_1 \overline{h}_1 \ot f_2\overline{h}_2 ) = 
   \frac{1}{|G|} \sum_{g\in G} 
 \ ^g(\ta\alpha)(f_1 \ot {}^{h_1}f_2)\,\overline{h_1h_2} 
\end{equation}
(see Remark~\ref{conversion}).
If $\mu_1$ integrates,
then the above formula
gives the coefficient of $t$ in the product 
of $f_1\overline{h}_1, f_2\overline{h}_2$ in the corresponding
deformation of $S(V)\# G$ over $\CC[t]$.
In particular, when $\alpha$ is constant, the formula defines 
the first multiplication map $\mu_1$ of a deformation
of $S(V)\#G$ arising from a graded Hecke algebra.  
(We see directly in that case that $\mu_1$ must lower degree by $2$.)

%%%%%%%%%%%%%%%%%%%%%%%%%%%%%%%%%%%%%%%%%%%%%%%%%%%%%%%%%%%%

\subsection*{Graded Hecke algebras to Deformations
and Hochschild Cocycles}
Consider a graded Hecke algebra defined by $\kappa$.
Define $\alpha$ in $C^2\cong\Hom_{\CC}(\Wedge^2 V, \CC G)$ 
by $$\alpha(v\wedge w)= \kappa(v,w) \in \CC G\ .$$
Theorem~\ref{gamma-prime} implies that 
$\alpha$ defines a constant Hochschild cocycle
in $\HH^2(S(V)\#G)$,
since $\kappa$ defines the first multiplication map $\mu_1$
of a deformation 
(i.e., a cocycle in $\HH^2(S(V)\# G)$)
by Theorem~\ref{degree-mu-i}.
(Alternatively, we may use Theorem~\ref{constantcocycles} above.)
The second author~\cite{Witherspoon2}
showed how to define functions $\mu_i$ ($i\geq 1$)
giving the corresponding deformation of $S(V)\# G$ over $\CC[t]$.
But the construction of the $\mu_i$ is iterative, involving repeated
applications of the relations in the graded Hecke algebra.
Our closed formula~(\ref{productoftwo}) improves this description by
giving the multiplication map $\mu_1$ 
in terms of quantum differentiation.

%%%%%%%%%%%%%%%%%%%%%%%%%%%%%%%%%%%%%%%%%%%%%%%%%%%%%%%%%%%%

%%%%%%%%%%%%%%%%%%%%%%%%%%%5
\subsection*{Faithful versus nonfaithful actions}

We end by pointing out that it is not sufficient merely to consider $G$ modulo
the kernel of its representation in this theory:
The Hochschild cohomology of $S(V)\#G$ for
$G$ acting nonfaithfully on $V$ requires extra care.
As an example, we explicitly point out 
the contribution from the kernel of the representation of $G$ on $V$
to the space of graded Hecke algebras.
(The effect of the kernel on the ring structure of cohomology
under cup product
is described in~\cite{paper2}.) 
Theorem~\ref{constantcocycles} and Remark~\ref{degree-two}
imply (also see~\cite[Theorem~1.9]{RamShepler} 
and~\cite[Corollary~8.17]{SheplerWitherspoon}):
\vspace{2ex}
\begin{cor}
The parameter space of graded Hecke algebras is isomorphic to
$$
\begin{aligned}
\bigoplus_{\substack{g\in\mathcal C\\ \rule[-.5ex]{0ex}{2ex}\codim V^g=2\\ 
\det h|_{(V^g)^{\perp}} =1,\ \forall h\in Z(g)}} 
\hspace{-3ex}\bigl(\CC\,\overline{g}\ot\vol_g^{\perp}\bigr)
\ \ \ \oplus \ 
\displaystyle\bigoplus_{\substack{g\in\mathcal C\\ V^g=V\rule[-.5ex]{0ex}{2ex}}} 
    \bigr(\,\overline{g}\ot \Wedge^{2}V^*\bigl)^{Z(g)} \ ,
\end{aligned}
$$
where $\vol_g^{\perp}$ is any fixed choice of volume form
on $((V^g)^{\perp})^*$, and $\mathcal C$ is a set of representatives
of conjugacy classes of $G$.
\end{cor}
\vspace{1ex}
\vspace{1ex}

%%%%%%%%%%%%%%%%%%%%%%%%%%%%%%%%%%%%%%%%%%%%%%%%%%%%%%%%%%%%%%%%%%%%%55
%%%%%%%%%%%%%

%%%%%%%%%%%%%%%%%%%%%%%%%%%%%%%%%%%%%%%%%%%%%%%%%%%%%%%%%%%%%%5
%%%%%%%%%%%%%%%%%%%%%%%%%%%%%%%%%%%%%%%%%%%%%%%%%%%%%%%%%%%%%%%%%%%%%
%%%%%%%%%%%%%%%%%%%%%%%%%%%%%%%%%%%%%%%%%%%%%%%%%%%%%%%%%%%%%%%%%%%%%


\begin{thebibliography}{99999}


\bibitem{BGS} A.\ Beilinson, V.\ Ginzburg, and W.\ Soergel,
``Koszul duality patterns in represenation theory,'' 
J.\ Amer.\ Math.\ Soc.\ 9 (1996), no.\ 2, 473--527.

\bibitem{BlockGetzler} J.\ Block and E.\ Getzler, ``Quantization of foliations,''
Proceedings of the XXth International Conference on Differential Geometric Methods
in Theoretical Physics, Vol.\ 1, 2 (New York, 1991), 471-487, World Sci.\ Publ.,
River Edge, NJ, 1992.

\bibitem{BravermanGaitsgory} A.\ Braverman and D.\ Gaitsgory, 
``Poincar\'e-Birkhoff-Witt theorem for quadratic algebras of Koszul type,''
J.\ Algebra 181 (1996), no.\ 2, 315--328.

\bibitem{CaldararuGiaquintoWitherspoon} A.\ C\u{a}ld\u{a}raru, A.\ Giaquinto,
and S.\ Witherspoon, ``Algebraic deformations arising from orbifolds
with discrete torsion,'' J.\ Pure Appl.\ Algebra 187 (2004), no.\ 1--3, 51--70.

\bibitem{Cherednik} I.\ Cherednik, ``Double affine Hecke algebras and
Macdonald's Conjectures," Ann.\ of Math.\ (2) 141 (1995), no.\ 1, 191--216.

\bibitem{Drinfeld} V.\ G.\ Drinfeld, ``Degenerate affine Hecke algebras
and Yangians,'' Funct.\ Anal.\ Appl.\ 20 (1986), 58--60.

\bibitem{EtingofGinzburg} P.\ Etingof and V.\ Ginzburg, ``Symplectic
reflection algebras, Calogero-Moser space, and deformed Harish-Chandra
homomorphism,'' Invent.\ Math.\ 147 (2002), no.\ 2, 243--348.

\bibitem{Farinati} M.\ Farinati, ``Hochschild duality, localization,
and smash products,'' J.\ Algebra 284 (2005), no.\ 1, 415--434.

\bibitem{Gerstenhaber} M.\ Gerstenhaber, ``On the deformation of rings
and algebras,'' Ann.\ Math.\ 79 (1964), 59--103.

\bibitem{GinzburgKaledin} V.\ Ginzburg and D.\ Kaledin, ``Poisson
deformations of symplectic quotient singularities,'' Adv.\ Math.\ 186 (2004),
no.\ 1, 1--57.

\bibitem{Gordon} I.\ Gordon ``On the quotient ring by diagonal invariants,"
Invent.\ Math.\ 153 (2003), no.~3, 503--518.

\bibitem{HalboutTang} G.\ Halbout and X.\ Tang, ``Noncommutative Poisson
structures on orbifolds,'' Trans.\ Amer.\ Math.\ Soc.\ 362 (2010),
no.\ 5, 2249--2277.

\bibitem{Kontsevich} M.\ Kontsevich, ``Deformation quantization of 
Poisson manifolds,'' Lett.\ Math.\ Phys.\ 66 (2003), no.\ 3, 157--216.

\bibitem{Loday} J.-L.\ Loday, {\em Cyclic Homology}, 2nd.\ ed., Springer-Verlag,
Berlin, 1998.

\bibitem{Lusztig3} G.\ Lusztig, ``Affine Hecke algebras and their graded
version,'' J.\ Amer.\ Math.\ Soc.\ 2 (1989), no.\ 3, 599--635.

\bibitem{RamShepler} A.\ Ram and A.V.\ Shepler, ``Classification
of graded Hecke algebras for complex reflection groups,'' Comment.\ Math.\
Helv.\ 78 (2003), 308--334.

\bibitem{SheplerWitherspoon} A.V.\ Shepler and S.\  Witherspoon,
``Hochschild cohomology and graded Hecke algebras,'' Trans.\ Amer.\ 
Math.\ Soc.\ 360 (2008), no.\ 8, 3975--4005.

\bibitem{paper1} A.V.\ Shepler and S.\ Witherspoon,
``Quantum differentiation and chain maps of bimodule complexes,''
to appear in Algebra and Number Theory.

\bibitem{paper2} A.V.\ Shepler and S.\ Witherspoon,
``Finite groups acting linearly: Hochschild cohomology and the cup product,''
to appear in Adv.\ Math.

\bibitem{Stefan} D.\ \c{S}tefan, ``Hochschild cohomology on Hopf Galois
extensions,'' J.\ Pure Appl.\ Algebra 103 (1995), 221--233.

\bibitem{Weibel} C.\ A.\ Weibel, {\em An Introduction to Homological Algebra},
Cambridge Studies in Adv.\ Math.\ 38, Cambridge Univ.\ Press, Cambridge, 1994.

\bibitem{Witherspoon2} S.\ Witherspoon, ``Twisted graded Hecke algebras,''
J.\ Algebra 317 (2007), 30--42.

\bibitem{Xu} P.\ Xu, ``Noncommutative Poisson algebras,'' Amer.\ J.\ Math.\ 116 (1994),
no.\ 1, 101--125.

\end{thebibliography}
\end{document}